\documentclass[11pt,reqno,english]{amsart}

\usepackage{amsmath}
\usepackage{amssymb}
\usepackage{amscd}
\usepackage{amsthm}
\usepackage{verbatim}
\usepackage{latexsym}

\usepackage{verbatim}
\usepackage[T2A]{fontenc}
\usepackage[cp1251]{inputenc}
\usepackage{babel}

\newcommand{\att}[2]{\genfrac{}{}{0pt}{}{#1}{#2}}

\newcommand{\cal}{\mathcal}

\newtheorem{th}{Theorem}
\newtheorem{prop}{Proposition}
\newtheorem{cor}{Corollary}

\newtheorem{lem}{Lemma}
\newtheorem{conj}{Conjecture}

\theoremstyle{definition}
\newtheorem{rem}{Remark}

\begin{document}
 \title{
 Random linear combinations of functions from~$L_1$}
 \author{P.G. Grigoriev}\thanks{This work was supported by grants
``Leading Scientific Schools of the Russian Federation'' 00-15-96047
and RFFI 02-01-00315.\\
\hspace*{\parindent}
\!\!\!\! This paper is to be published in Mathematical Notes (\it Matematicheskie zametki).}
 \maketitle
\date{}
\begin{abstract}
In the paper I study properties of random polynomials
with respect to a general system of functions.
Some lower bounds for the mathematical expectation
of the uniform and recently introduced integral-uniform norms
of random polynomials are established.
\end{abstract}

\section{Introduction}

In the article we research random polynomials
(linear combinations of functions with random coefficients)
of the type
\begin{equation}\label{polinom2}
F_n(\omega,x):=\sum_{i=1}^n \xi_i(\omega)f_i(x),
\end{equation}
where $\{f_i\}_{1}^n$ is a system of functions on a measure
space~$(X,\mu)$ and $\{\xi_i\}_{i}^n$ is a set of independent
random variables on a probability space~$(\Omega,{\sf P})$
(``random coefficients''). Our interest will be concentrated on
the properties of the random variable
$\|F_n(\omega,\cdot)\|_{B(X)}$, where ${\|\cdot\|_B}$ is a norm in
a certain space~$B$ consisting of functions on~$(X,\mu)$.

Note, when ${B=L_\infty}$,
the formulated above research topic becomes the classical
problem of estimating
the expectation of the supremum of random polynomial~(\ref{polinom2}).
In~1954
Salem and Zygmund~\cite{sz}
established that
\begin{equation}\label{salz}
{\sf E}\big\|\sum_{k=-n}^nr_k(\omega)e^{ikx}\big\|_\infty
\equiv \frac1{2^{2n+1}}
\hspace{-.6em}
\sum_{\att{\text{over all choices}}{\text{of signs:}\
\varepsilon_k=\pm1}}
\hspace{-.4em}
\Big\|\sum_{k=-n}^n\varepsilon_ke^{ikx} \Big\|_\infty\!
\asymp\sqrt{n\log n},
\end{equation}
where $r_k(\omega)=\text{sign}\sin(2^{k+1}\pi\omega)$
are the Rademacher functions on~$[0,1]$ and
${\sf E}$ denotes the mathematical expectation.
However, by the virtue of Khinchin's inequality it follows that
$$
{\sf E}\big\|
\sum_{k=-n}^n  r_k(\omega)e^{ikx}\big\|_p
\asymp_p\sqrt{n},
\quad\ 1\le p<\infty.$$
These facts demonstrate that
the embeddings of the spaces of trigonometric polynomials
into $L_\infty$ and~$L_p$ with~${p<\infty}$
are qualitatively different.

By now a number of estimates for the uniform norm
of random polynomials~(\ref{polinom2})
with various constraints on~$\{\xi_i\}_1^n$ and~$\{f_i\}_{1}^n$
have been established.
Note while the upper estimates have
many important applications in analysis and probability,
the applications of the lower estimates are rather rare.
Nevertheless, proofs of the lower estimates
are more challenging and usually need involving
some specific properties of systems~$\{f_i\}$
and~$\{\xi_i\}$.
E.g. the proof of the lower bound in~(\ref{salz})
in~\cite{sz} essentially relied on the facts that the trigonometric
functions is a system of {\it characters}
(i.e. formulae like
${\cos(\alpha+\beta)}={\cos\alpha\cdot\cos\beta
-\sin\alpha\cdot\sin\beta}$
were used)
and that for distribution of linear combinations of Rademacher functions
an exponential estimate holds.
In~\cite{ks}
one can find a short proof of a sharp lower estimate
for the random polynomial~$\sum_{k}\xi_k\cos kx$, where
$\{\xi_k\}$ are independent Gaussian variables,
this proof relies on the fact that the matrix
$\big(\cos \frac{2\pi}n jk\big)_{j,k=0}^{n-1}$
is a multiple of an orthogonal matrix and that
uncorrelated Gaussian variables are independent.
Proofs of the upper estimates normally use a
simple averaging argument and
need essentially weaker constraints on
system~$\{f_i\}$.
Usually, if an upper estimate for the expectation
of a norm of the random polynomial~(\ref{polinom2})
takes place, then essentially the same estimate
holds if we substitute in~(\ref{polinom2}) the functions~$f_i$ with
their absolute values~$|f_i|$.
In this paper, continuing the work~\cite{ja3},
there is established
a lower estimate for the uniform norm of
random polynomials with respect to an abstract
functional system~$\{f_i\}_1^n$ from~$L_1$,
provided that the system satisfies a weak condition.

In monographs~\cite{ledouxtal}, \cite{marcuspis}
one can find a deep theory which enables to estimate
the expectation of the uniform norm of random polynomials~(\ref{polinom2}),
provided $\{f_i\}$  is a system of characters of a
locally compact Abelian group restricted to a compact
neighborhood~$V$ of the group identity element.
The basic method for estimating the supremum
of a random process in~\cite{ledouxtal},~\cite{marcuspis}
(in particular, of a random polynomial~(\ref{polinom2}))
is reduction of the problem to
a problem of estimating the $\varepsilon$-entropy of $V$
with respect to a metric induced by the random process.
This method was originally introduced by Dudley~\cite{dud} and
Sudakov~\cite{sud1},~\cite{sud2}
and later was developed by Fernique, Marcus, Pisier, Talagrand
and others (see \cite{ledouxtal},
\cite{marcuspis}, \cite{sud3},~\cite{fern}).
To obtain a lower estimate for the maximum
of a random process using this method, one has to apply a variant
of Slepian's lemma~\cite{sl} for Gaussian vectors.
This lemma enables one to estimate from below the probability
$$
{\sf P}\big\{\max_{1\le k\le m}
\sum_{i=1}^n\xi_if_i(x_k) \ge\alpha\big\},
\qquad\{x_k\}_1^m \text{\ is a net in\ }X,
$$
provided that {\bf all} normalized inner products
(cosines of angles):
$$
\frac{(W_{x_j},W_{x_k})}{|W_{x_j}|\cdot|W_{x_k}|}
,\quad
{j\ne k},\quad~
\text{for vectors}\quad
W_{x_j}:=\big(f_i(x_j)\big)_{i=1}^n\in{\Bbb C}^n,
\quad~ j=1,\dots,m
$$
are sufficiently small and~$\{\xi_i\}$
are independent Gaussian variables.
Demonstrating the existence of such a net~$\{x_j\}$
is usually a separate non-trivial problem.
Moreover, since Slepian's lemma can be applied
only to Gaussian vectors
there arise some serious difficulties
with transfer of the estimates for
random polynomials with Gaussian coefficients
to the case of non-Gaussian~$\{\xi_i\}$.

In 1995 Kashin and Tzafriri \cite{kt1}, \cite{kt2}, \cite{kt3}
introduced another method for obtaining lower estimates
of the uniform norm of random polynomials.
In particular, in~\cite{kt1} it was shown that the lower estimate
in the Salem-Zygmund theorem~(\ref{salz})
stays true for random polynomials with respect
to an {\bf arbitrary} orthonormal system,
provided that the system is uniformly bounded in~$L_3$.
This approach relies on a version of the central limit
theorem with an estimate of the error term.
To apply this method
it is not necessary to estimate
{\bf all} angles between the vectors~$W_{x_j}$,
instead of that it suffices to demonstrate
that cosines of these angles  are {\bf small on average}.

In \cite{kt2}
Kashin and Tzafriri defined the following norm
\begin{equation}\label{norma}
 \|f\|_{m,\infty}:=\int_X...\int_X \max\{|f(x_1)|,...,|f(x_m)|\}
d\mu(x_1)...d\mu(x_m),
\end{equation}
where $f$ is an integrable function on a measure space $(X,\mu)$,
${\mu(X)=1}$.
As in author's works~\cite{ja2}, \cite{ja3}
let us call it {\it integral-uniform norm}
for this norm.
While obtaining the lower estimate
for the expectation of the uniform norm
of the random polynomial~(\ref{polinom2})
in~\cite{kt1},~\cite{kt3},
a similar estimate  for the
{\it integral-uniform} norm
with parameter ${m\asymp n^{1/2+\varepsilon}}$
was {\it de~facto} proved.
It is worth mentioning that estimates for
the norm $\|f\|_{m,\infty}$
(or rather for the {\it family} of norms with~${m\ge1}$)
are self-interesting,  since the values of the norms
carry quite full information about the distribution of function~$f$:
$\lambda_f(t)={\mu\{x\in X:|f(x)|\ge t\}}$.
Indeed, it is easy to see that
\begin{equation}\label{norma2}
\|f\|_{m,\infty}=\int_0^\infty
\big(1-(1-\lambda_f(t))^m\big)dt
\end{equation}
and
${\|f\|_1\equiv\|f\|_{1,\infty}}$\,${\le\|f\|_{m,\infty}\le\|f\|_\infty}$,
moreover,
${\|f\|_{m,\infty}\!\to\|f\|_\infty}$ as ${m\to\infty}$.
In Section~\ref{auxiliary} (Th.~\ref{kee}) it will be shown that
$$
\|f\|_{m,\infty}
\asymp
\sup_{\att{\Delta\subset X}{\mu\Delta=1/m}}
\Big\{m\int_\Delta|f|\,d\mu\Big\}
$$
provided there exists a subset of $X$ of measure $1/m$.

Using a simple modification of
the method from~\cite{kt1},~\cite{kt3},
the author~\cite{ja3}
proved the following

\smallskip
\noindent
{\bf Theorem A.}\label{genth}
{\it
Let
 $\{f_i\}_{i=1}^n$ and $\{\xi_i\}_{i=1}^n$
be systems of functions defined on probability
spaces  $(X,\mu)$ and $(\Omega,{\sf P})$  respectively,
satisfying:
\begin{itemize}
\item[\bf (a)]  $\|f_i\|_2=1$ and $\|f_i\|_{2+\varepsilon}\le M$
for all
$i=1,\ldots,n$
with some constants~${M>0}$, ${\varepsilon>0}$.
\item[\bf(b)]\label{bstar}
$
\big\|\sum\limits_{i=1}^n c_i f_i\big\|_2\le
M R^p
\big(\sum\limits_{i=1}^n|c_i|^2\big)^{1/2}$
for all sets of coefficients~$\{c_i\}_{1}^n$
with some constants $p\in[0,\frac12)$, ${M>0}$, where
$
R\equiv
R(\{a_i\}_1^n):=\frac{(\sum_{i=1}^n|a_i|^2)^2}{\sum_{i=1}^n|a_i|^4}$
and
$\{a_i\}_{1}^n$ is a {\bf fixed}
set of complex coefficients.
\item[\bf (c)]
$\{\xi_i\}_1^n$
is a system of independent variables such that
${\sf E}\xi_i=0$, ${{\sf E}|\xi_i|^2=1}$ and
${({\sf E}|\xi_i|^{2+\varepsilon})\le M}$
\end{itemize}
Then there exist positive constants ${q=q(p)}$,
${C_j=C_j(p,M,\varepsilon)}$,
${j=1,2,3}$,
such that\footnote{\label{powersnosk}%
It is possible to show that in~(\ref{genth1})
one can choose the power parameter~$q$ arbitrarily from the interval
$\big(0,{(1-2p)/2}\big)$,
in this case the constants~${C_1}$, $C_2$
may depend on the choice of~$q$.
One can find a detailed proof of~(\ref{genth1})
with ${q=(1-2p)/4}$ in~\cite{jadis}.
}\begin{equation}\label{genth1}
{\sf P}\Big\{\omega\in\Omega :
\big\|\sum_{i=1}^n a_i\xi_i (\omega )f_i\big\|_{m,\infty}
\le C_1 \Big(\sum_{i=1}^n|a_i|^2\log P\Big)^{1/2}
\Big\}
\le \frac{C_2}{P^q},
\end{equation}
where $P=\min(m,R)+1$ and
$R=R(\{a_i\})$ are defined in condition~{\bf(b)}.
It obviously implies
\begin{equation}\label{genth2}
{\sf E}\Big\|\sum_{i=1}^n a_i\xi_i f_i\Big\|_{m,\infty}
\ge C_3\Big(\sum_{i=1}^n|a_i|^2\log P\Big)^{1/2}.
\end{equation}
(In the case of small values of $m$ or $R$ the inequality (\ref{genth2})
easily follows from Khinchin's and H\"older's inequalities and
the trivial estimate:
 $\|\cdot\|_{m,\infty}\ge \|\cdot\|_1$.)
}
\pagebreak[1]

Theorem~A
is a simple generalization of the results from \cite{kt1}--\cite{kt3},
where the \mbox{$\|\cdot\|_\infty$-norm} of random
polynomials was estimated.
There estimates of type (\ref{genth1}), (\ref{genth2})
for the uniform norm were proved in the case~${p=0}$
and was noted that if~${a_i\equiv1}$,
these estimates could be generalized for the case~${p\in[0,1/2)}$.
It was also noted that the established estimates stay true
for the integral-uniform norm
with parameter ${m\asymp R^{1/2}}$.

Roughly speaking, Theorem~A implies that
an estimate of type~(\ref{genth2})
for a random polynomial~(\ref{polinom2})
does {\bf not} take place
only if the system~$\{f_i\}$
(of normalized in~$L_2$ functions)
is ``significantly far'' from an
orthonormal one (e.g. if the functions of the system
converge to a fixed function).

Note when
${a_i\equiv1}$ (${R(\{a_i\})=n}$),
applying Lemma~\ref{geomlem} (see below),
one can easily deduce the condition~{\bf(b)}
for functions~$\{f_i\}$, satisfying ${\|f_i\|_2=1}$,
from the condition:
\begin{itemize}\label{bshtrih}
\item[\bf(b$'$)]  $\|\sum_{i=1}^n \varepsilon_if_i\|_2\le
Mn^{\frac12+\beta}$
for all choices of signs $\varepsilon_i=\pm1$\\
with some ${\beta\in[0,\frac12)}$, ${M>0}$.
\end{itemize}

\smallskip
In~\cite{ja3}
it was also shown that provided ${m\le R}$
and some additional constraints on~$\{\xi_k\}_1^n$
the lower estimate~(\ref{genth2})
is sharp in sense of order.
It was demonstrated by the following

\noindent
{\bf Theorem B.}\label{upper}
{\it
Let $\{\xi_k\}_1^n$ be a system of independent variables
such that the following exponential estimate holds:
\begin{equation*}\label{exp}
{\sf P}\Big\{\big|\sum_{k=1}^n c_k\xi_k\big|> t\big(\sum_{k=1}^n
c_k^2\big)^{1/2}\Big\} \le C_4e^{-t^2C_5}
\end{equation*}
for all sets of coefficients $\{c_k\}_{1}^n$
with some constants~$C_4$, ${C_5>0}$.
Then there exists a constant $C_6=C_6(C_4,C_5)>0$
such that
\begin{equation*}\label{upper1}
{\sf E}\big\|\sum_{k=1}^n\xi_k f_k\big\|_{m,\infty} \le
C_6\big\| \big(\sum_{k=1}^n|f_k|^2\big)^{1/2}\big\|_{m,\infty}
\cdot\sqrt{1+\log m}
\end{equation*}
{\bf for all}
systems of functions ${\{f_k\}_1^n\subset L_1(X,\mu)}$
(${\mu X=1}$)
and for all~${m\ge1}$.
}
\pagebreak[1]

\smallskip
Note that the  condition imposed in Theorem~A,
which requires the functions~$\{f_i\}$
to be {\it uniformly bounded} in~$L_{2+\varepsilon}$,
looks rather unnatural since the theorem
provides {\it lower} estimates.
Roughly speaking, the reason for necessity
of such a condition is that the uniform boundness
of system~$\{f_i\}$ in~$L_{2+\varepsilon}$,
along with its orthonormality
(or a weaker condition~{\bf(b$'$)}
combined with~${\|f_i\|_2=1}$),
ensures that the ``essential'' supports of functions~$f_i$
mutually intersect ``strongly enough''.
To convince that the condition of uniform boundness
cannot be simply omitted from Theorem~A
consider the example of the functional system ${f_i:=\sqrt{n}\chi_i}$
on $[0,1]$, where $\chi_i$ are the indicators of
the intervals $\big(\frac{i-1}n,\frac in\big)$.

The main target of the article is to generalize Theorem~A
in the particular case~${a_i\equiv1}$
for the random polynomials of type~(\ref{polinom2})
with respect to a system of functions~$\{f_i\}_1^n$, ${\|f_i\|_1=1}$,
{\it which are not necessarily bounded} in~$L_p$, ${p>1}$.
To avoid extreme functional systems, such as one in the
previous paragraph, the conditions~{\bf(a)}
and {\bf(b$'$)}
are substituted for
\begin{itemize}\label{l1d}
\item[\bf (d)]  $\|f_i\|_1=1$  for all $i=1,\cdots,n$ and \\
$\|\sum_{i=1}^n\theta_if_i\|_1\le
Mn^{\frac12 +p}$
for all choices of signs~$\{\theta_i\}_1^n$, $\theta_i=\pm1$
\end{itemize}
with some constants $p\in[0,\frac1{12})$, $M>0$.
The main result is the following
\pagebreak[1]

\begin{th}\label{l1th}
Let $\{f_i\}_{i=1}^n$  be a system of functions
on a probability space $(X,\mu)$
which satisfies the condition~{\bf(d)} with~$p\in[0,\frac1{12})$.
Let~$\{\xi_i\}_{i=1}^n$ be independent variables defined on another
probability space~$(\Omega,{\sf P})$, satisfying
${{\sf E}\xi_i=0}$, ${{\sf E}|\xi_i|^2=1}$ and
${{\sf E}|\xi_i|^3\le M^3}$.
Then there exist some constants
$q'=q'(p)>0$, $C'_j=C'_j(p,M)>0$, $j=1,2,3$,
such that whenever ${m\le n}$
\begin{equation}\label{l1th1}
{\sf P}\Big\{\omega\in\Omega :
\| F_n(\omega,x)
\|_{m,\infty} \le C'_1
\sqrt{n\cdot(1+\log m)} \Big\} \le
{C'_2}{m^{-q'}}
\end{equation}
and, consequently,
\vspace{0em plus 0em minus .2em}
\begin{equation}\label{l1th2}
{\sf E}\|\sum_{i=1}^n \xi_i f_i\|_{m,\infty} \ge C'_3
\big(n\cdot(1+\log m)\big)^{1/2},
\end{equation}
where $F_n$ is a random polynomial
defined by~(\ref{polinom2}).
(For small $m$ inequality~(\ref{l1th})
follows from Khinchin's inequality.)
\end{th}

This result provides a new estimate not only for the integral-uniform
norm but also for the uniform norm of random
polynomials~(\ref{polinom2}):

\begin{cor}\label{l1uniform}
Let functions~$\{f_i\}_1^n$
and random variables~$\{\xi_i\}_1^n$ satisfy
the conditions of Theorem~\ref{l1th}.  Then for the {\em uniform}
norm of the random polynomial~(\ref{polinom2}) the following
estimate holds
\begin{equation*}
{\sf P}\Big\{\omega\in\Omega :
\| F_n(\omega,x)\|_{\infty} \le C'_1
\sqrt{n\cdot(1+\log n)} \Big\} \le
{C'_2}{n^{-q'}}.
\end{equation*}
\end{cor}
\smallskip
\pagebreak[1]

Note that Theorem~\ref{l1th}, being stronger than
Theorem~A in some environments, is weaker than that
in the two aspects:
first, it cannot be applied to polynomials
of type $\sum_k a_k\xi_k(\omega) f_k(x)$ with
an arbitrary choice of non-random coefficients~$\{a_i\}$;
second, it imposes the condition~{\bf(d)}
with parameter~$p<\frac1{12}$, while in the corresponding
condition for Theorem~A
(condition~{\bf(b)} provided ${a_i\equiv1}$ or~{\bf(b$'$)})
it suffices to have~${p<1/2}$.
I think that the constraint~${p<\frac1{12}}$
can be relaxed\footnote{\samepage%
The condition~{\bf(d)} with parameter~${p=1/2}$
is obviously satisfied for {\bf all} functional systems~$\{f_i\}_1^n$,
${\|f_i\|_1=1}$.
In particular, it holds for the trivial system~${f_i\equiv1}$
for which by Khinchin's inequality the
$L_\infty$-norm of random polynomial with respect to that system
is of order~$\sqrt n$.
Thus,
it makes sense to think about possible validity of estimates
of type~(\ref{l1th1}),~(\ref{l1th2}) only in the case~${p<1/2}$.
}
(see Conjecture~\ref{conja1} below).
In~\cite{jadis}
a conjecture about possible
generalization of Theorem~\ref{l1th} for the
case of polynomials~$\sum_k a_k\xi_k f_k$
(with non-trivial $a_k$) is formulated.
\pagebreak[1]

The organization of the paper is as follows.
In Section~\ref{auxiliary} we shall prove several auxiliary results
which were explicitly or implicitly used in~\cite{kt1},
\cite{kt3},~\cite{ja3}.
We gather these results with the intention of making the method
developed in~\cite{kt1}, \cite{kt3},~\cite{ja3}
easier to understand, apply and make appropriate alterations,
e.g. for generalizing of Theorem~A
for the case of specific values of the power parameter~$q$
(see Footnote~\ref{powersnosk} and Remark~\ref{powerq} below).
In Section~\ref{secl1}
we shall prove Theorem~\ref{l1th}.
In Section~\ref{applic}
Theorem~\ref{l1th} will be applied to give
a partial solution for
a Functional Analysis problem formulated by
Montgomery-Smith and Semenov in~\cite{sem}
and to estimate the Marcinkiewicz norm of random polynomials.
In Section~\ref{applic} we shall also formulate
two hypotheses,  concerning the
potential generalizations of Theorem~\ref{l1th}.

I would like to express my special gratitude to B.S.~Kashin,
whose advisement led me to establish the results of the paper,
I thank also E.M.~Semenov and A.M.~Zubkov for
valuable remarks and discussions.
%
%
%
%

\section{Auxiliary Results}\label{auxiliary}

{\bf Integral-Uniform Norm.}
Let us check that the definitions~(\ref{norma}) and~(\ref{norma2})
of the integral-uniform norm are identical.
It is well-known that ${\|g\|_1=\int_0^\infty\lambda_g(t)\,dt}$
for every function ${g\in L_1(Y,\nu)}$,
where $\lambda_g(t):={\nu\{y:|g(y)|\ge t\}}$ is
the distribution of~$g$.
Thus, to prove equivalence of~(\ref{norma}) and~(\ref{norma2})
it suffices to notice that the function
$$
g(\bar x)=\max\big\{|f(x_1)|,...,|f(x_m)|\big\},
\ \quad \bar x=(x_j)_1^m\in X^m=:Y,
$$
has the distribution: $\lambda_g(t)=1-(1-\lambda_f(t))^m$.

By inequality ${\max(|a|,|b|)\le|a|+|b|}$
it is easy to see that
\begin{equation}\label{qw4}
\|f\|_{m,\infty}\le m\|f\|_1\quad\mbox{and}\quad
\|f\|_{n,\infty}\le\frac{n+1}m\|f\|_{m,\infty}
\ \quad\mbox{for all $f\in L_1(X)$, $m\le n$}.
\end{equation}
For the indicator~$\chi_\Delta$ of a set~$\Delta\subset X$
the identity~(\ref{norma2}) implies
$
\|\chi_\Delta\|_{m,\infty}={1-(1-|\Delta|)^m}
$
(here and further we denote  $|\Delta|\equiv\mu\Delta$).
Thus, if~$m\ge c|\Delta|^{-1}$, then $\|\chi_\Delta\|_{m,\infty}\ge C(c)$
with a constant {$C(c)>0$}.

For an integrable function~$f\in L_1$
define the following ``relative'' norms:
\begin{align}\label{norma3}
\|f\|_m^* &:=
\sup_{\Delta\subset X}
\Big\{\frac{1-(1-|\Delta|)^m}{|\Delta|}
\int_\Delta|f|d\mu\Big\};\\
\label{norma4}
\|f\|'_m &:=
\sup_{\att{\Delta\subset X}{\mu\Delta=1/m}}
\Big\{m\int_\Delta|f|d\mu\Big\}.
\end{align}
Right side of (\ref{norma4}) is well-defined only
if there exists an event $e\subset X$
such that $\mu (e)=1/m$.
However, we will use this norm only if
$X=[0,1]$ with standard Lebesque measure.

The ${\|\cdot\|^*_m}$- and ${\|\cdot\|'_m}$-norms are equivalent to ${\|\cdot\|_{m,\infty}}$-norm,
this is proved in the following result, which
was implicitly established by the author in~\cite{ja3}.

\begin{th}\label{kee}
For all functions~$f\in L_1[0,1]$
the following inequalities take place
\begin{equation}\label{kl}
(1-\frac1e)\|f\|'_m \le
\|f\|_{m}^* \le
\|f\|_{m,\infty} \le
2\|f\|'_m.
\end{equation}
\end{th}

\smallskip
\noindent
{\bf Proof.}
Let us show that $2\|f\|'_m\ge\|f\|_{m,\infty}$.
Let
$\Delta^*\subset [0,1]$ satisfy $|\Delta^*|=1/m$ and
$$m\int_{\Delta^*}|f|=\|f\|'_m.$$
(It is easy to see that such a~$\Delta^*$ exists
though not necessarily unique).
We have
$$\|f\|_{m,\infty}\le
\|f\cdot\chi_{\Delta^*}\|_{m,\infty}+
\|f\cdot(1-\chi_{\Delta^*})\|_{m,\infty}.$$
By (\ref{qw4}) we estimate
$\|f\cdot\chi_{\Delta^*}\|_{m,\infty}
\le m\int_{\Delta^*}|f|=\|f\|'_m$.
From the extremality of $\Delta^*$
we get
$$
\|f\cdot(1-\chi_{\Delta^*})\|_{m,\infty}
\le\|f\cdot(1-\chi_{\Delta^*})\|_\infty
\le\frac1{|\Delta^*|}\int_{\Delta^*}|f|=\|f\|'_m.
$$
Thus, the inequality $\|f\|_{m,\infty}\le2\|f\|'_m$ is proved.

Let us check now that for every $\Delta\subset [0,1]$
the following inequality holds:
\begin{equation}\label{s11}
\|f\|_{m,\infty}\ge
\frac{1-(1-|\Delta|)^m}{|\Delta|}\int_\Delta|f|d\mu.
\end{equation}
Obviously, it suffices to consider the case when $f$
vanishes outside~$\Delta$
(${\text{supp}(f)\subset\Delta}$).
By~(\ref{norma2}) we get
\begin{align*}
\|f\|_{k+1,\infty}-\|f\|_{k,\infty}
&=\int_0^\infty \!\big(1-(1-\lambda_f(t))^{k+1}\big)dt-
\int_0^\infty \!\big(1-(1-\lambda_f(t))^k\big)dt \\
&=\int_0^\infty \lambda_f(t)\big(1-\lambda_f(t)\big)^k dt \\
&\ge (1-|\Delta|)^k \int_0^\infty \!\lambda_f(t) dt=(1-|\Delta|)^k \|f\|_1.
\end{align*}
Summing these inequalities from~$k=1$ to~$k=m-1$, we get
$$
\|f\|_{m,\infty}-\|f\|_1
\equiv
\|f\|_{m,\infty}-\|f\|_{1,\infty}
\ge\|f\|_1\sum_{k=1}^{m-1}(1-|\Delta|)^k,$$
which implies
\vspace{0pt plus 0pt minus .2em}
$$
\|f\|_{m,\infty}\ge\|f\|_1\frac{1-(1-|\Delta|)^m}{|\Delta|}.$$
Now, to prove~(\ref{s11}) it remains to notice that
$\|f\|_1=\int_\Delta |f|$
(recall $\text{supp}f\subset\Delta$).
Therefore, $\|f\|_{m,\infty}\ge\|f\|_m^*$.

The inequality
${\|f\|^*_m\ge(1-e^{-1})\|f\|'_m}$
obviously follows from the fact that ${(1-\frac1m)^m<e^{-1}}$.
The proof of Theorem~\ref{kee} is completed.
\pagebreak[2]
\medskip

{\bf A Fact from Geometry.}\nopagebreak

\begin{lem}\label{geomlem}
(See Lemma~1 in \cite{ja3}  or Lemma~2.1 in \cite{jadis}).
Let $\{w_i\}_{i=1}^n$
be a set of vectors in a linear space equipped with
a norm~$\|\cdot\|$ (or a semi-norm), satisfying $\|w_i\|=1$ and
\begin{equation}\label{geom}
\|\sum_{i=1}^n \theta_i w_i\|\le c\cdot n^{\frac12+\beta}
\nopagebreak
\end{equation}
for all choices of signs $\{\theta_i\}_{i=1}^n$, $\theta_i=\pm1$,
with some constants
${\beta\in[0,1/2)}$, ${c>0}$.
Then
$$
\big\|\sum_{i=1}^n a_i w_i\big\|\le
C(c)\, n^{\frac14+\frac\beta2}
\big(\sum_{i=1}^n a_i^2\big)^{\frac12}
\eqno(\ref{geom}')
$$
for all sets of coefficients  $\{a_i\}_1^n$.
This estimate is sharp, i.e. there exist
vectors~$\{w_i\}_1^n$,
a norm ${\|\cdot\|}$, and coefficients~$\{a_i\}_1^n$
such that~(\ref{geom}) takes place, but
estimate~(\ref{geom}\,$'$)  is sharp in sense of order.
\end{lem}

Geometrically Lemma~\ref{geomlem} claims that
the convex hull of the set $B^d_\infty\cup (n^{1/2+\beta}\cdot B^d_1)$
has the inscribed sphere with radius of order~$n^{1/4+\beta/2}$,
here $B^d_\infty$ denotes $d$-\nolinebreak{dimensional} cube
whose vertices have coordinates~$\pm1$ and
${B^d_1:=\{(v_k)\in{\Bbb R}^d:\sum_1^d|v_k|\le1\}}$.
\pagebreak[1]

\smallskip
{\bf Transfer Lemmas.}
To prove the main result we need the following
lemmas which enable us to transfer one property of
the multidimensional  normal distribution
(Lemma~\ref{lem2})
to the case of an abstract multidimensional distribution.

\begin{lem}\label{lem1}
Let~$\Omega_j\subset\Omega$, $j=1,\dots,m$,
be events, satisfying
\begin{equation*}
(1-\kappa)
\sum_{j,k=1}^m{\sf P}\big({\Omega_j\cap \Omega_k}\big)
\le
\Big(\sum_{j=1}^m {\sf P}(\Omega_j)\Big)^{2}
\end{equation*}
with some~$\kappa\in(0,1)$. Then \,
${\sf P}\big(\bigcup_{j=1}^m\Omega_j\big)
\ge1-\kappa.$
\end{lem}
\pagebreak[1]

\vspace{.0em plus .3em minus .1em}
\begin{lem}\label{lem2}
Let~$\{h_j\}_{j=1}^m$ be a set of Gaussian variables
such that ${{\sf E}h_j=0}$,
\hspace{0em plus .3em minus .0em}
${ {\sf E}h_j^2=D_j\ge r^2>0}$,
\hspace{0pt plus .3em minus .0em}
${{\sf E}h_jh_k=v_{j,k}}$
and each pare
$(h_j,h_k)$, ${j\ne k}$
has a 2-dimensional normal distribution with density
\vspace{0pt plus 0pt minus .2em}
$$
\phi_{0,V_{j,k}}(Y):=\frac1{2\pi(\det V_{j,k})^{1/2}}
\exp\big\{-\frac12(Y,V_{j,k}^{-1}Y)\big\},\qquad\ Y=(y_1,y_2),
$$
$$
\text{where}\qquad\
V_{j,k}=\left(
\begin{array}{cc}
D_j & v_{j,k}\\
v_{j,k} & D_k
\end{array}
\right)
\quad\text{is a covariance matrix.}
$$
Assume also that there exist some constants
$R\ge1$, $c_0$,~${\delta>0}$ such that
$$
\frac1{m^2}
\sum_{\att{j,k=1}{j\ne k}}^m |v_{j,k}|\le c_0R^{-\delta}r^2.
$$
Then for arbitrary choice of ${\alpha<\delta^{1/2}}$
there exists an {\em independent of~$R$} constant
${C_7=C_7(c_0,\alpha,\delta)}$ such that
\begin{equation}\label{lem22}
\sum_{{j,k=1}}^m
{\sf P}\big(\Psi_j\cap\Psi_k\big)
\le(1+C_7P^{-q_0}) \Big(\sum_{j=1}^m {\sf P}(\Psi_j)\Big)^2,
\end{equation}
where
${P:=\min(R,m)+1}$,
$
\Psi_j=\Psi_j(\alpha):={\big\{h_j>\alpha\sqrt{D_j\log P}\big\}}$
and
$q_0=q_0(\alpha,\delta)=
{\min\big(\frac13(\delta-\alpha^2),\frac32\alpha^2\big)}>0$.
\end{lem}
\pagebreak[1]


\vspace{.0em plus .3em minus .01em}
\begin{lem}\label{lem3tex}
Let~$(\eta_j)_1^m$ and $(h_j)_1^m$ be random vectors with
identical first and second moments:
$$
{\sf E}\eta_j  ={\sf E}h_j=0,\quad\
{\sf E}\eta_j\eta_k={\sf E}h_jh_k=v_{j,k},\quad\
D_j\equiv v_{j,j}\ge r^2>0.
$$
Moreover, let~$(h_j)_1^m$ be
a Gaussian vector whose covariance
matrix satisfies the assumption of Lemma~\ref{lem2}
with parameters~$r$, $R$,~$c_0$,~${\delta>0}$.
Assume also that
there exist some positive constants
${\delta_i\le1}$, $M_i$ ${i=1,2,3}$,
and $\alpha<\alpha_0:=
\big(\min(\delta,\delta_1,\delta_2,\delta_3)\big)^{1/2}$
such that
\begin{align}\label{lem3usl1}
\big|{\sf P}\big(U_j(\alpha)\big)-
{\sf P}\big(\Psi_j(\alpha)\big)\big|
&\le M_1P^{-\delta_1},\quad \ 1\le j\le n; \\
\label{lem3usl2}
\big|{\sf P}\big(U_j(\alpha)\cap U_k(\alpha)\big)-
{\sf P}\big(\Psi_j(\alpha)\cap\Psi_k(\alpha)\big)\big|
&\le M_2P^{-\delta_2},\quad \ (j,k)\in\sigma,
\end{align}
where $P=\min(R,m)+1$ and
\begin{align*}
\Psi_j=\Psi(\alpha)&:=
 \big\{\omega:h_j(\omega)>\alpha\sqrt{D_j\log P}\big\};\\
U_j=U_j(\alpha)&:=\big\{\omega:\eta_j(\omega)>\alpha\sqrt{D_j\log P}\big\};
\end{align*}
and the index set~$\sigma\subset\{(j,k):j\ne k,1\le j,k\le m\}$
satisfies ${|\sigma|\ge m^2(1-M_3P^{-\delta_3})}$.\\
\hspace*{\parindent plus 0em minus .1em}
Then
there exists a constant
${C_7'=C_7'(\alpha,c_0,\delta,\delta_i,M_i)}$, ${i=1,2,3}$,
such that
\begin{equation}\label{lem33tex}
\sum_{j,k=1}^m
{\sf P}\big(U_j(\alpha)\cap U_k(\alpha)\big)
\le(1+C_7'P^{-q'})
\Big(\sum_{j=1}^m {\sf P}\big(U_j(\alpha)\big)\Big)^2,
\end{equation}
where the power parameter
$q'=
\min\big\{
\frac13(\delta-\alpha^2),\frac32\alpha^2,
\frac56\min\limits_{i=1,2,3}(\delta_i-\alpha^2)\big\}>0$.
\end{lem}
\pagebreak[1]

\noindent
{\bf Remark on Lemmas \ref{lem1}--\ref{lem3tex}.}
Lemma~\ref{lem1}, being a generalization
of the Borel-Cantelli Lemma, is a known and important result.
One can find similar statements, e.g. in
\cite{cher},
\cite{sprind}, \cite{kt1}.
Lemmas~\ref{lem2} and \ref{lem3tex},
as far as I am aware of, are formulated for the first time,
however the ideas for their proofs
have been thoroughly borrowed
from~\cite{kt1},~\cite{kt3}.
Lemma~\ref{lem3tex}
serves to transfer the estimates of type~(\ref{lem22})
to the ``non-Gaussian case.''
Note that the main results of~\cite{kt1}, \cite{kt3},~\cite{ja3}
could be easier proved and perceived with the help of
Lemmas~\ref{lem1} and~\ref{lem3tex}.

Note also that from Lemmas~\ref{lem1} and~\ref{lem2}
one can easily derive for Gaussian random variables~$\{h_i\}_1^m$,
satisfying the assumption of Lemma~\ref{lem2},
the following estimate:
$$
{\sf P}(\max\limits_{1\le j\le m}h_j  > \alpha r\sqrt{\log m})
\ge{1-Cm^{-q_0}} \qquad \text{for}\quad m\le R,
$$
where $R$ and $r$ are from the statement of Lemma~\ref{lem2}
and $\alpha$, $q_0$,~$C$ are some positive constants.
This fact links Lemma~\ref{lem2} with
the results of Slepian~\cite{sl}
and \v{S}idak~\cite{sidak1}, \cite{sidak2},
devoted to estimating the
distribution of maximum of Gaussian vectors with
a non-trivial covariance matrix.
\pagebreak[1]

\noindent
{\bf Proof of Lemma~\ref{lem1}.}
Let $\chi_j$ be  the
indicators of the events ${\Omega_j\subset\Omega}$
and let ${\zeta:=\sum_{j=1}^n\chi_j}$. Then
${\sf E}|\zeta|=\sum_{j=1}^m{\sf P}(\Omega_j)$,
$
{\sf E}|\zeta|^2=\sum_{j,k=1}^m
{\sf P}\big({\Omega_j\cap \Omega_k}\big)
$
and \,${\text{supp}\,\zeta}=\bigcup_{j=1}^m\Omega_j$.
Applying the Cauchy-Schwarz inequality, we get
$$
(1-\kappa)^{1/2}\big({\sf E}|\zeta|^2\big)^{1/2}\le
{\sf E}|\zeta|\le\big({\sf E}|\zeta|^2\big)^{1/2}
\Big({\sf P}\Big(\bigcup_{j=1}^m\Omega_j\Big)\Big)^{1/2}.
$$
This proves Lemma~\ref{lem1}.
\pagebreak[1]

\smallskip
\noindent
{\bf Proof of Lemma~\ref{lem2}.}
Since the random variables~$h_j$
are normal we have
\begin{equation}\label{lem2nomer1}
{\sf P}(\Psi_j)=\frac1{\sqrt{2\pi D_j}}
\int_{\alpha\sqrt{D_j\log P}}^\infty
e^{-{y^2}/({2D_j})}\,dy=
\frac1{\sqrt{2\pi}}\int_{\alpha\sqrt{\log P}}^\infty
e^{-{y^2}/2}\,dy.
\end{equation}
Taking into account that
$
\int_z^\infty e^{-t^2/2}\,dt\asymp z^{-1}
e^{-z^2/2}$
when $z\ge1$,
we have
$C_\alpha^{-1} (\sqrt{\log P})^{-1}P^{-\alpha^2/2} \le
{\sf P}(\Psi_j)\le
C_\alpha (\sqrt{\log P})^{-1}P^{-\alpha^2/2}$
with a constant $C_\alpha$ depending only on $\alpha$. Thus
$$
\sum_{j=1}^m{\sf P}(\Psi_j)
\asymp_\alpha \frac m{\alpha\sqrt{\log P}}P^{-\alpha^2/2}.\eqno(*)$$
Thus, to prove~(\ref{lem22})
it suffices to check that
$$
\sum_{\att{j,k=1}{j\ne k}}^m
{\sf P}\big(\Psi_j\cap\Psi_k\big)
\le(1+C'_7P^{-q_0}) \Big(\sum_{j=1}^m {\sf P}(\Psi_j)\Big)^2.
\eqno(\ref{lem22}')
$$

Define the following index set:
$$
\sigma_1:=\big\{(j,k): 1\le j\ne k\le m,
|v_{j,k}|<\frac1{32}r^2\big\}. $$
Chebyshev's inequality for the set
\hspace{0em plus .2em minus .1em}
${\sigma_1^c=\{(j,k)\!:1\le j\ne k\le m\}\setminus\sigma_1}$
implies
$$
|\sigma_1^c|\le \frac{32}{r^2}\sum_{j\ne k}|v_{j,k}|
\le{32}c_0m^2R^{-\delta}.$$
Therefore,
$$
\sum_{(j,k)\in\sigma_1^c}{\sf P}(\Psi_j\cap\Psi_k)
\le {32}c_0m^2R^{-\delta}\max_{1\le j\le m}{\sf P}(\Psi_j).
$$
Thus, taking into account (\ref{lem2nomer1}), we get
$$
\sum_{(j,k)\in\sigma_1^c}{\sf P}(\Psi_j\cap\Psi_k)
\le 32c_0m^2R^{-\delta}\frac{\text{\rm const}}
{\alpha\sqrt{\log P}}P^{-\alpha^2/2}.
$$
Applying (*), we conclude that
$$\sum_{(j,k)\in\sigma_1^c}
{\sf P}\big(\Psi_j\cap\Psi_k\big)
\le K_2R^{-\delta+\alpha^2}
\cdot
\Big(\sum_{j=1}^m {\sf P}(\Psi_j)\Big)^2
\eqno(**)$$
with a constant $K_2(\alpha,c_0)>0$.
Thus, to prove~$(\ref{lem22}')$
we can neglect the summation over~$\sigma_1^c$ on the left-hand side.

Now, let us estimate the sum
$$
\Sigma_1:=\hspace{-.2em plus 1pt minus .7em}
\sum_{s=(j,k)\in\sigma_1}\hspace{-.2em plus 1pt minus .7em}
{\sf P}(\Psi_j\cap\Psi_k)
=\hspace{-.3em plus .1em minus 1.2em}
\int\limits_{\alpha\sqrt{D_j\log P}}^\infty
\int\limits_{\alpha\sqrt{D_k\log P}}^\infty
\sum_{s\in\sigma_1}
\frac{\exp\big\{\!-\frac12(Y,V_s^{-1}Y)\big\}}{2\pi\sqrt{\det V_s}}
dy_1dy_2.
$$
Changing the integration variables
$t_1=\frac{y_1}{\sqrt{D_j}}$,
$t_2=\frac{y_2}{\sqrt{D_k}}$,
we get
$$\Sigma_1=
\int_{\alpha\sqrt{\log P}}^\infty
\int_{\alpha\sqrt{\log P}}^\infty
\sum_{s\in\sigma_1}
\frac{\sqrt{D_jD_k}}{2\pi\sqrt{\det V_s}}
e^{-\frac12Q(t_1,t_2)}
dt_1dt_2,
$$
where $Q(t_1,t_2):=(Y,V_s^{-1}Y)$
is a quadratic form.
Evaluating the determinant,
we get ${\det V_s=D_jD_k-v_s^2}$ and
$$ V_s^{-1}= \frac1{D_jD_k-v_s^2}
\left(
\begin{array}{ccc}
D_k       &-{v_s}\\
 -{v_s}   &{D_j}
\end{array} \right).
$$
Thus, the coefficients of the quadratic form
\hspace{0pt plus .4em minus .1em}
${Q=a_st_1^2+a_st_2^2-2b_st_1t_2}$
\hspace{0pt plus .4em minus .1em}
are defined by
$$
a_s= \frac{D_jD_k}{D_jD_k-v_s^2};
\qquad\qquad
b_s= \frac{v_s\sqrt{D_jD_k}}{D_jD_k-v_s^2}.
$$
Notice also that $|v_s|\le r^2/32$, provided~$s\in\sigma_1$.
Therefore, taking into account that
$D_j\ge r^2$ for all~$j=1,\dots,m$, we have
$$a_s-b_s=\frac{\sqrt{D_jD_k}}{\sqrt{D_jD_k}+v_s}
\ge\frac12.
$$
Thus, we can estimate the values of the quadratic form~$Q$ as follows:
$$
Q(t_1,t_2)=a_st_1^2+a_st_2^2-2b_st_1t_2\ge
(a_s-b_s)(t_1^2+t_2^2)
\ge\frac12(t_1^2+t_2^2).
$$

Now, for arbitrary~${L>1}$ we can estimate
\begin{align*}
J(L)&:=\frac1{2\pi}
\int_{L\alpha\sqrt{\log P}}^\infty
\int_{\alpha\sqrt{\log P}}^\infty
\sum_{s\in\sigma_1}
\Big( \frac{D_jD_k}{D_jD_k-v_s^2}\Big)^{\frac12}
e^{-\frac12Q(t_1,t_2)}
dt_1dt_2\\
&\le
\int_{L\alpha\sqrt{\log P}}^\infty
\int_{\alpha\sqrt{\log P}}^\infty
\sum_{s\in\sigma_1}
e^{-\frac1{4}(t_1^2+t_2^2)}dt_1dt_2\\
&\asymp_\alpha
\frac{|\sigma_1|}{L\log P}e^{-\frac{\alpha^2}{4}(L^2+1)\log P}
\le
\frac{K_3m^2}{L\log P}P^{-\frac{\alpha^2(L^2+1)}{4}},
\end{align*}
where $K_3=K_3(\alpha)$ is a constant.
Choose $L=3$ and take into account~(*) to get
$$J(3)<K_4 P^{-\frac32{\alpha^2}}
\Big(\sum_{j=1}^m{\sf P}(\Psi_j)\Big)^2\eqno(***)$$
with a constant $K_4=K_4(\alpha)>0$.

Notice that
$$
\Sigma_1\le 2 J(3)+
\frac1{2\pi}
\hspace{-.1em plus .2em minus .2em}
\int_{\alpha\sqrt{\log P}}^{3\alpha\sqrt{\log P}}
\int_{\alpha\sqrt{\log P}}^{3\alpha\sqrt{\log P}}
\sum_{s\in\sigma_1}
\sqrt{a_s} \cdot
e^{-\frac12Q(t_1,t_2)}
dt_1dt_2.
$$
Now, in order to finish the proof of~(\ref{lem22})
it remains to compare the expression
$$
A:=\frac1{2\pi} \sum_{s\in\sigma_1}
\sqrt{a_s}
e^{-\frac12Q(t_1,t_2)}
$$
with the expression
$B:=(2\pi)^{-1}|\sigma_1|\exp\{-\frac12(t_1^2+t_2^2)\}$
on the square
${\alpha\sqrt{\log P}\le t_1,t_2\le 3\alpha\sqrt{\log P}}$.
In fact, if we show that
$A\le B(1+K_5 P^{-q_0})$
with a constant
${K_5(\alpha,\delta)>0}$,
then integrating this inequality we get
\begin{align*}
\frac1{2\pi}
& \int_{\alpha\sqrt{\log P}}^{3\alpha\sqrt{\log P}}
\int_{\alpha\sqrt{\log P}}^{3\alpha\sqrt{\log P}}
\sum_{s\in\sigma_1}
\sqrt{a_s} \exp\big\{-\frac12Q(t_1,t_2)\big\}dt_1dt_2<\\
&<(1+{K_5}{P^{-q_0}})
\int_{\alpha\sqrt{\log P}}^\infty
\int_{\alpha\sqrt{\log P}}^\infty
\frac{|\sigma_1|}{2\pi}
\exp\big\{-\frac12(t_1^2+t_2^2)\big\} dt_1dt_2\\
&= (1+{K_5}{P^{-q_0}})
\sum_{s\in\sigma_1}{\sf P}(\Psi_j)\cdot{\sf P}(\Psi_k).
\end{align*}
This inequality combined with (**) and (***)
would imply~(\ref{lem22})
and, thus, prove the lemma.

Split the index set $\sigma_1$ into the subsets
$$\sigma_i:=
\big\{s\in\sigma_1:2^{-i}r^2\le|v_s|<2^{-i+1}r^2\big\},
\qquad  \ \ i=6,7,\dots.$$
Clearly, $\sigma_1=\bigcup_{i\ge6}\sigma_i$.
Chebyshev's inequality for~$|\sigma_i|$ implies
$$
|\sigma_i|2^{-i}r^2\le\sum_{s\in\sigma_1}|v_s|
\le c_0r^2m^2R^{-\delta}
$$
\vspace{0pt plus 0pt minus .2em}%
and, consequently,
\vspace{0pt plus 0pt minus .1em}
$$|\sigma_i|\le\min\big(m^2,
2^ic_0m^2 R^{-\delta}\big).$$
Taking into account that $D_j\ge r^2$ for $s\in\sigma_i$ we have
\begin{multline*}
\big|(a_s-1)(t_1^2+t_2^2)-2b_st_1t_2\big|
\le\big(|a_s-1|+|b_s|\big)(t_1^2+t_2^2)=\\
=\frac{|v_s|}{\sqrt{D_jD_k}-|v_s|}(t_1^2+t_2^2)
<2^{2-i}(t_1^2+t_2^2)\le2^{7-i}\alpha^2\log P
\end{multline*}
in the domain
$\alpha\sqrt{\log P}\le t_1,t_2\le 3\alpha\sqrt{\log P}$.
Moreover, it is easy to see that if $s\in\sigma_i$ and $i\ge6$, then
$$\sqrt{a_s}=
\Big(
\frac{{D_jD_k}}{D_jD_k-v_s^2}
\Big)^{1/2}
=\Big(1+
\frac{v_s^2}{D_jD_k-v_s^2}
\Big)^{1/2}
\le1+2^{3-2i}.$$

Gathering all these facts, we get
\vspace{0em plus 0pt minus .2em}
$$\quad A\le\frac{B}{|\sigma_1|}S,$$
$$
\text{where}\quad\
\hspace{0em plus 1em minus .4em}
S:=
\hspace{-.1em plus 1em minus .3em}
\sum_{i=6}^\infty |\sigma_i|(1+2^{3-2i})
e^{\frac12 2^{7-i}\alpha^2\log P}
=
\hspace{-.1em plus 1em minus .4em}
\sum_{i=6}^{[2q_0\log_2 P]}
\hspace{-.1em plus 1em minus .5em}
+
\hspace{-.2em plus 1em minus .5em}
\sum_{i>[2q_0\log_2 P]}^\infty
\hspace{-.3em plus 1em minus .5em}=:S_1+S_2.
$$
Taking into account that
$q_0=\min\big(\frac13(\delta-\alpha^2),\frac32\alpha^2\big)$,
we estimate
$$ S_1 \le \hspace{-.3em plus 1em minus .6em}
\sum_{i=6}^{[2q_0\log_2 P]}
\hspace{-.3em plus 1em minus .6em}
2^ic_0m^2R^{-\delta}\!\cdot
(1+2^{-9})P^{\alpha^2}
\hspace{-.1em plus 1em minus .2em}
\le 4c_0 m^2 P^{\alpha^2+2q_0}R^{-\delta}
\hspace{-.1em plus 1em minus .2em}
\le K_6 R^{-q_0}|\sigma_1|,
$$
where $K_6=K_6(c_0)>0$ is a constant.
\pagebreak[1]

To estimate~$S_2$ notice that
$P^{2^6\alpha^2P^{-2q_0}}<(1+K_{7}P^{-q_0})$,
where  $K_{7}(\alpha,q_0)>0$ is a constant. Thus,
$$
S_2\le
\sum_{i>[2q_0\log_2 P]}
|\sigma_i|(1+8P^{-4q_0})(1+K_7P^{-q_0})\le
|\sigma_1|(1+K_8 P^{-q_0}),$$
where   $K_{8}=8+9K_7$. Therefore,
$$
S\le\big(1+(K_6+K_8)P^{-q_0}\big)|\sigma_1|$$
and
$A\le B(1+K_5P^{-q_0})$. This completes
the proof of the inequality~(\ref{lem22}) and Lemma~\ref{lem2}.
\pagebreak[1]

\begin{rem}\label{remlem2}
In Lemma~\ref{lem2} the
power parameter~$q_0(\alpha,\delta)$
is not optimally chosen.
It is not difficult to show that the inequality~(\ref{lem22})
stays true, provided
$q_0={q_0(\alpha,\delta,\varepsilon)=
\min\big\{\frac12(\delta-\alpha^2),\frac32\alpha^2\big\}
-\varepsilon>0}$ for arbitrarily small $\varepsilon$.
(In order to check this it suffices to draw a sharper estimate
for the sum~$S$).
However, in this case the constant~$C_7$
in the inequality~(\ref{lem22}) would depend also on~${\varepsilon>0}$.
Moreover, somewhat more advanced modification of the proof
enables to derive~(\ref{lem22}) with parameter
${q_0=\frac12(\delta-\alpha^2)-\varepsilon}$.
\end{rem}
\pagebreak[1]

\smallskip
\noindent
{\bf Proof of Lemma \ref{lem3tex}.}
Since the random vector~$(h_j)_1^m$ satisfies the assumption
of Lemma~\ref{lem2} it follows that the estimate~(\ref{lem22})
for the events~$\{\Psi_j\}_{j=1}^m$
with the power parameter ${q_0=\frac13(\delta-\alpha^2)}$ holds.
Thus, to prove~(\ref{lem33tex})
it suffices to show that for arbitrary~$\alpha<\alpha_0$
the following inequalities take place:
\begin{align}
\label{lem3t1}
\sum_{j=1}^m{\sf P}\big(\Psi_j(\alpha)\big)
& \le
 \big(1+L_1P^{-q_1}\big)
 \sum_{j=1}^m{\sf P}\big(U_j(\alpha)\big);\\
\label{lem3t2}
 \sum_{j,k=1}^m\!{\sf P}\big(U_j(\alpha)\cap U_k(\alpha)\big)
& \!\le
\!\sum_{j,k=1}^m\!
{\sf P}\big(\Psi_j(\alpha)\cap \Psi_k(\alpha)\big)+
\frac{L_2}{P^{q_2}}
 \Big(\sum_{j=1}^m\!{\sf P}\big(\Psi_j(\alpha)\big)\Big)^2,
\end{align}
where the constants~$q_1$,~$q_2\ge q'$
and~$L_1$,~${L_2>0}$ depend only on~$\alpha$, $c_0$,
$\delta$, $\delta_i$, $M_i$, ${i=1,2,3}$.

Notice that
\begin{equation}\label{psi4}
{\sf P}\big(\Psi_j(\alpha)\big)=
\frac1{2\pi D_j}
\int_{\alpha\sqrt{D_j\log P}}^\infty
e^{-\frac{y^2}{2D_j}}dy
\asymp_\alpha \frac{P^{-\alpha^2/2}}{\alpha\sqrt{\log P}}.
\end{equation}
Thus, for all~${\alpha^2<\delta_1}$ we have
(see~(\ref{lem3usl1}))
\begin{align}
\label{lem3t3}
{\sf P}\big(U_j(\alpha)\big)
&\asymp_\alpha
{\sf P}\big(\Psi_j(\alpha)\big);\\
\notag
M_1P^{-\delta_1}
&\le L_3 P^{-\delta_1+\alpha^2}
\cdot {\sf P}\big(\Psi_j(\alpha)\big)
\end{align}
with a constant~$L_3(\alpha,\delta_1,M_1)$.
Taking into account~(\ref{lem3usl1}),
we get the inequality~(\ref{lem3t1})
with parameter~${q_1=\delta_1-\alpha^2}$.

In order to prove~(\ref{lem3t2}) let us recall that according
to the Lemma assumption
$$
|\sigma^c|\equiv
\text{card}\big\{(j,k):j\ne k,(j,k)\notin\sigma\big\}
\le M_3m^2P^{-\delta_3}.
$$
Thus, for all~$\alpha^2<\delta_1$ we get
(see~(\ref{lem3t3}))
$$
\sum_{(j,k)\in\sigma^c}{\sf P}(U_j\cap U_k)
\le M_3m^2P^{-\delta_3}\max_{1\le j\le m}{\sf P}(U_j)
\le L_4 m^2P^{-\delta_3}
\frac{P^{-\alpha^2/2}}{\sqrt{\log P}}
$$
with a constant~$L_4(\alpha,\delta_1,\delta_3,M_1,M_3)$.
Using~(\ref{psi4}) and (\ref{lem3t3}), for
${\alpha^2<\min\{\delta_1,\delta_3\}}$ we get
$$
\sum_{(j,k)\in\sigma^c}
{\sf P}\big(U_j(\alpha)\cap U_k(\alpha)\big)
\le L_5P^{-\delta_3+\alpha^2}
\Big( \sum_{j=1}^m{\sf P}\big(\Psi_j(\alpha)\big)\Big)^2.
$$
Hence, to prove~(\ref{lem3t2})
we can neglect the summation over~$\sigma^c$
on the left-hand side.

Notice, when~${\alpha^2<\delta_1}$ the estimate~(\ref{lem3t3})
implies
$$
\sum_{j=1}^m{\sf P}\big(U_j(\alpha)\big)
\le L_5 P^{-1+\alpha^2}
\Big(\sum_{j=1}^m{\sf P}\big(\Psi_j(\alpha)\big)\Big)^2
$$
with a constant $L_5(\alpha,\delta_1,M_1)$.
Therefore, to prove~(\ref{lem3t2})
we can also neglect  summation over the
pares~${\big\{(j,k):j=k\big\}}$
on the left-hand side.

Now, to prove~(\ref{lem3t2}) it remains to notice
that the assumption~(\ref{lem3usl2}) implies
$$
\sum_{(j,k)\in\sigma}
{\sf P}\big(U_j(\alpha)\cap U_k(\alpha)\big)
\le\sum_{(j,k)\in\sigma}
{\sf P}\big(\Psi_j(\alpha)\cap\Psi_k(\alpha)\big)
+M_2m^2P^{-\delta_2}.
$$
Taking into account (\ref{psi4}), for~$\alpha^2<\min\{\delta_1,\delta_2\}$
we can estimate the error term as follows:
$$
 M_2m^2P^{-\delta_2}\le
 L_6 P^{-\frac56({\delta_2-\alpha^2}) }
\Big(\sum_{j=1}^m{\sf P}\big(\Psi_j(\alpha)\big)\Big)^2,
$$
where $L_6=L_6(\alpha,\delta_1,\delta_2,M_1,M_2)$ is a constant.
Thus, the inequality~(\ref{lem3t2}) is proved with the constants
${L_2=L_4+L_5+L_6}$  and
${q_2=\min\{\frac56({\delta_2-\alpha^2}),\delta_3-\alpha^2\}}$.

The inequalities~(\ref{lem3t1}),~(\ref{lem3t2}),
combined with~(\ref{lem22}), prove~(\ref{lem33tex})
with ${q'=\min\{q_0,q_1,q_2\}}$,
where ${q_0=\min(\frac13(\delta-2\alpha^2),\frac32\alpha^2)}$
is the power parameter in~(\ref{lem22}).
The proof of Lemma~\ref{lem3tex} is completed.

\begin{rem}\label{powerq}
A simple modification to the proof of Theorem~A in~\cite{ja3},
which would involve Lemmas~\ref{lem1}--\ref{lem3tex},
enables one to establish the estimate~(\ref{genth1})
with the power parameter~$q$ arbitrarily chosen
from~${\big(0,\frac3{11}({1-2p})\big)}$.
Moreover, if we took unto account Remark~\ref{remlem2}
and made some simple refinements in Lemmas~\ref{lem2} and~\ref{lem3tex},
then  we could prove~(\ref{genth1})
with the parameter~$q$
from the interval ${\big(0,({1-2p})/2\big)}$,
however in this case the constants~$C_1$, $C_2$
would depend on~$q$.
In~\cite{jadis} one can find a detailed proof of Theorem~A
for the case ${q=(1-2p)/4}$.
\end{rem}
\pagebreak[1]

{\bf Central Limit Theorem.}
For the proof of the main result we need to apply
a version of 2-dimensional central limit theorem with
an estimate of the error term.
We shall use in one- and two-dimensional case the
following result due to Rotar'~\cite{rotar}
(or see Corollary~17.2 in~\cite{bhat}):
\begin{prop}\label{rot}
Let
$\{X_i\}_{i=1}^N$
be a set of independent random vectors
in ${\Bbb R}^d$, satisfying ${\sf E}X_i=0$, $1\le i\le N$,
then
$$\sup_{A\in {\mathcal C}} |P_N(A)-\Phi_{0,V}(A)|\le K_1(d)
N^{-1/2}m_3\lambda^{-3/2},$$
where $P_N(A)$
is the probability of the event that
$N^{-1/2}\sum_{i=1}^NX_i$ belongs to the set~$A$,
${\mathcal C}$~denotes the class of all Borel convex sets in ${\Bbb R}^d$,
$K_1(d)<\infty$ is a constant,
$$m_3:=\frac1N\sum_{i=1}^N{\sf E} |X_i|^3,$$
$\lambda$ is the smallest eigenvalue of the matrix
$V=N^{-1}\sum_{i=1}^N \text{\rm cov} (X_i)$,
$\text{\rm cov} (X_i)$ denotes the covariance matrix of vector~$X_i$,
finally, $\Phi_{0,V}$ denotes the Gaussian measure on~${\Bbb R}^d$
with the density
$$\phi_{0,V}(Y):=(2\pi)^{-d/2}(\det V)^{-1/2}
\exp\big\{-\frac12(Y,V^{-1}Y)\big\},\qquad \ \ Y\in {\Bbb R}^d.$$
\end{prop}
\pagebreak[1]

\section{Proof of Theorem~\ref{l1th}}\label{secl1}

Theorem~\ref{l1th} is a direct corollary of
more general Theorem~\ref{l1th}$'$
and Khinchin's inequality.

\smallskip
\noindent
{\bf Theorem {\ref{l1th}}$'$.}
{\it
Let $\{f_i\}_{i=1}^n$
be a system of functions on a probability space~$(X,\mu)$, satisfying
\begin{itemize}
\item[\bf({d}$'$)]
$\|f_i\|_1=1$  for all $i=1,\cdots,n$;\\
$\big\|\sum_{i=1}^n \theta_i f_i\big\|_1\le Mn^{\frac 12+p_1}$
for all choices of signs $\theta_i=\pm1$;\\
$\big\|(\sum_{i=1}^n|f_i|^2)^{1/2}\big\|_1\le Mn^{\frac12+p_2}$,
where $M$, ${p_1, p_2\ge0}$ are  some constants, satisfying
 $p_1+2p_2<\frac12$ and~${p_2<\frac1{12}}$.
\end{itemize}
Let $\{\xi_i\}_{i=1}^n$
be a system of independent random variables on another
probability space~$(\Omega,{\sf P})$, satisfying
${{\sf E}\xi_i=0}$, ${{\sf E}|\xi_i|^2=1}$ and
${{\sf E}|\xi_i|^3\le M^3}$.
Then whenever~${m\le n}$ for random polynomial~(\ref{polinom2})
the estimates~(\ref{l1th1}) and~(\ref{l1th2}) hold
with some constants
$q=q(p_1,p_2)>0$, $C_j=C_j(p_1,p_2,M)>0$, ${j=1,2,3}$.
(These constants, of course, are not the same as the
constants from the statement of Theorem~A.)
}
\medskip

To deduce Theorem~\ref{l1th} from Theorem~\ref{l1th}$'$
it suffices to notice that validity of condition~{\bf(d)},
combined with integrated over~${x\in X}$
Khinchin's inequality (e.g. see \cite{ks})
for the sum ${\sum r_i(\omega) f_i(x)}$
with fixed~$x$, where
$r_i$ are the Rademacher functions, implies
$$\big\|\Big(\sum_{i=1}^n |f_i|^2\Big)^{\frac12}\big\|_1
\le \text{const}\cdot{\sf E}\int_X
\big|\sum_{i=1}^n r_i(\omega) f_i(x)\big|d\mu(x)
\le M'(M)n^{\frac12+p}.
$$
Thus, if condition~{\bf({d})} holds with parameter $p=p_1$,
the condition~{\bf({d}$'$)} holds with $p_2\le p_1$
and, consequently, if~${p_1<\frac1{12}}$,
then  ${p_1+2p_2<1/2}$ automatically.

Theorem~\ref{l1th}$'$, like Theorem~A,
is based on the central limit theorem
(Proposition~\ref{rot}), however its proof requires
essentially subtler preparatory work.
Roughly speaking, the reason for this
is that while the basis functions~$\{f_i\}$
were supposed to be {\it uniformly} bounded in~$L_{2+\varepsilon}$
it was possible to find a ``sufficiently large'' set~${E\subset X}$
such that
\begin{equation}\label{tridva}
\frac{\sum_{i=1}^n|f_i(x)|^3}{\big(\sum_{i=1}^n|f_i(x)|^2\big)^{3/2}}
\le\text{\rm{const}}\cdot n^{-\epsilon},\qquad\ x\in E,
\end{equation}
with some $\epsilon>0$.
This inequality is needed to estimate the error term
after the application of the central limit theorem.
Such a trick was used to prove~(\ref{genth1}),~(\ref{genth2})
in~\cite{ja3}
(and before that it had been uses in~\cite{kt1},~\cite{kt3}
to prove similar estimates for the uniform norm).

However, in the assumptions of Theorem~\ref{l1th}$'$
there may be no point~${x\in X}$ such that (\ref{tridva}) holds.
As a corresponding  example take the
functions ${f_i:=r_i+n^{q}\chi_i}$, where
$r_i$ are the Rademacher functions,
$\chi_i$ are the characteristic functions of the intervals
$\big(\frac{i-1}n,\frac in\big)$ and~${\frac13<q<\frac7{12}}$.
(You can normalize $f_i$ in~$L_1$ to make the
example more appropriate for Theorem~\ref{l1th}$'$.)
For such functions the condition~{\bf(d)} holds
with~$p=\max\{0,q-\frac12\}<\frac1{12}$,
however the inequality~(\ref{tridva}) fails a.e. on~$[0,1]$.
Thus, one cannot directly apply Proposition~\ref{rot},
as well as {\it the other} versions of the central limit theorem,
for the sums $\sum_1^n\xi_i f_i(x)$.
Nevertheless, the condition~{\bf(d$'$)}
reserves a possibility for
``sufficiently large'' set of points $x\in X$
to pick out ``long enough'' subsum $\sum_{i\in{\cal I}_x}\xi_if_i(x)$
for which the estimate~(\ref{tridva}) holds and, consequently,
the central limit theorem can be applied.
It turns out that the index set~${\cal I}_x$ may depend on~$x$,
for this reason there arise some difficulties with  the transfer
of estimates for subsums to the case of original polynomial.
\pagebreak[1]


\smallskip
\noindent
{\bf Proof of Theorem \ref{l1th}$'$.}
Without loss of generality assume ${2\le m\le n^{\varepsilon_1}-1}$
with some fixed
${\varepsilon_1=\varepsilon_1(p_1,p_2)\in(0,\frac12)}$,
whose value will be explicitly set later.
Moreover, since ${\|f\|_{m,\infty}}$
depends only on the distribution of $f$ (see~(\ref{norma2}))
and $L_1[0,1]$ contains equimeasurable copy of vector $(f_j)_1^n$
we can suppose that $X=[0,1]$
with standard Lebesque measure.
The last means that we can use Theorem~\ref{kee}.

\smallskip
\noindent
{\bf Step 1.}
For each $x\in X$ define the index set
$${\cal A}(x)\equiv {\cal A}_x:=\Big\{k:1\le k\le n,
\; |f_k(x)|\le n^{-\frac12-\varepsilon_1}\sum_{j=1}^n|f_j(x)|\Big\}.$$
By Chebyshev's inequality we get
$|{\cal A}_x^c|\le n^{\frac12+\varepsilon_1}$
and $|{\cal A}_x|\ge {n-n^{\frac12+\varepsilon_1}}$\!.
For each $x\in X$ define a smaller index set:
${\Lambda(x)\equiv\Lambda_x\subset {\cal A}_x}$ such that
${|\Lambda_x|=n-[n^{\frac12+2\varepsilon_1}]}$
and at every~$x$ the set~$\Lambda_x$ indexes the
${n-[n^{\frac12+2\varepsilon_1}]}$ least values of~$|f_k(x)|$.
In order to define the set~$\Lambda_x$ formally and
provide it with an additional property that for all
$k\in\{1,\cdots,n\}$ the set ${\{x:k\in\Lambda(x)\}}$
is $\mu$-measurable we use the following inductive procedure.
Let
\begin{align*}
k_1(x) &:=
 \min\big\{ k\in\{1,\dots,n\}: |f_k(x)|\ge|f_i(x)|~
 \forall i=1,\dots,n\big\};\\
{\cal K}_1(x) &:=\big\{1,\dots,n\big\}\setminus\{k_1(x)\}.
\end{align*}
Assume that $k_l(x)$, ${\cal K}_l(x)$ for
$l=1,\dots,j-1$ are defined and set
\begin{align*}
k_j(x) &:=
 \min\big\{ k\in{\cal K}_{j-1}(x): |f_k(x)|\ge|f_i(x)|~
 \forall i\in{\cal K}_{j-1}(x)\big\};\\
{\cal K}_j(x) &:={\cal K}_{j-1}(x)\setminus\{k_j(x)\}.
\end{align*}
It is easy to show that the indices~$k_j(x)$, ${j=1,\dots,n}$,
are measurable functions of~$x$. Set
$
\Lambda_x :=\big\{1,\dots,n\big\}\setminus
{\cal K}_{[n^{1/2+2\varepsilon_1}]}(x).
$
Clearly,
$|\Lambda_x|={n-[n^{\frac12+2\varepsilon_1}]}$,
${\Lambda_x\subset{\cal A}_x}$ and
$$
\int_X\sum_{k\in \Lambda_x^c}|f_k(x)|d\mu(x) \le\int_X
\sqrt{|\Lambda_x^c|}\cdot
\Big(\sum_{k=1}^n|f_k|^2\Big)^{1/2}d\mu(x) \le
Mn^{\frac34+p_2+\varepsilon_1}.$$
and, consequently,
\vspace{0em plus 0em minus .2em}
$$
\int_X\sum_{k\in \Lambda_x}|f_k(x)|d\mu(x)\ge
n-Mn^{\frac34+p_2+\varepsilon_1}.
$$
Non-triviality of this estimate will be ensured by
the choice of~$\varepsilon_1$,
satisfying ${\frac34+p_2+\varepsilon_1<1}$.
\pagebreak[1]

\smallskip
\noindent
{\bf Step 2.}
Set
$$
~X':=\Big\{x\in X:
\sum_{k\in\Lambda_x}|f_k(x)|\ge \frac13\sum_{k=1}^n|f_k(x)|\Big\}$$
and notice that
$
\int_{X\setminus X'}\sum_{k\in\Lambda_x}|f_k|
\le\frac13\int_X\sum_{k=1}^n|f_k|
=n/3$.
Therefore, for sufficiently large\footnote{As usual,
\samepage
the case of small $n$ can be dealt by increasing or
reducing some constants.}
$n\ge n_0(\varepsilon_1+p_2,M)$ we have
$$
\int_{X'}
\sum_{\Lambda_x}|f_k|\ge n-Mn^{\frac34+p_2+\varepsilon_1}-\frac
n3\ge\frac n2.$$
H\"{o}lder's inequality for the function
$$F_2(x):=\Big(\sum_{k\in\Lambda_x}|f_k(x)|^2\Big)^{1/2}$$
implies
$$\int_{X'}F_2(x)\,d\mu(x)\ge
\int_{X'}|\Lambda_x|^{-\frac12}\cdot\sum_{k\in \Lambda_x}|f_k|\,d\mu(x)
\ge \frac{\sqrt n}2.$$

Notice, if $x\in X'$ and an index set~$I_x\subset {\cal A}_x$
satisfies~$|I_x|\ge|\Lambda_x|$, then by the
definitions of~$\Lambda_x$ and~${\cal A}_x$ we get
\begin{align*}
\sum_{k\in I_x}|f_k(x)|^3 &\le \Big(\sum_{k\in I_x}|f_k(x)|^2\Big)
\cdot\max_{k\in I_x}\big\{|f_k(x)|\big\}
\le \Big(\sum_{k\in I_x}|f_k(x)|^2\Big)
\cdot n^{-\frac12-\varepsilon_1}\sum_{k=1}^n|f_k(x)|\\
&\le n^{-\frac12-\varepsilon_1}\Big(\sum_{k\in I_x}|f_k(x)|^2\Big)
\cdot 3\Big(\sum_{k\in\Lambda_x}|f_k(x)|\Big)\\
&\le 3n^{-\frac12-\varepsilon_1}
|\Lambda_x|^{\frac12}\Big(\sum_{k\in I_x}|f_k(x)|^2\Big)^{3/2}
\le 3n^{-\varepsilon_1} \Big(\sum_{k\in I_x}|f_k(x)|^2\Big)^{3/2}.
\end{align*}
Thus
\begin{equation}\label{tridva2}
\frac{\sum_{k\in I_x}|f_k(x)|^3}
{\big(\sum_{k\in I_x}|f_k(x)|^2\big)^{3/2}}\le 3n^{-\varepsilon_1}.
\end{equation}
We shall need this inequality to estimate the
error term in the central limit theorem.
\pagebreak[1]

\smallskip
\noindent
{\bf Step 3.}
Define the sets
$$E_\ell:=\Big\{x\in X':\frac{\sqrt n}4 2^{\ell-1}
\le F_2(x)<\frac{\sqrt n}4 2^\ell \Big\},
\qquad\ \ell=1,2,\dots.$$
Assuming $n$ sufficiently large to ensure that
$$\frac{\sqrt n}2\le\int_{X'}F_2\le Mn^{\frac12+p_2},$$
for arbitrary $\varepsilon_2>0$ we have
\begin{equation}\label{st3}
\mu\Big\{x\in X':F_2(x)\ge\frac{\sqrt n}4 2^{n^{\varepsilon_2}} \Big\}
\le\frac4{\sqrt n} 2^{-n^{\varepsilon_2}}\int_{X'}F_2\,d\mu
<
4Mn^{p_2} 2^{-n^{\varepsilon_2}}.
\end{equation}
Notice
\vspace{0em plus 0.01em minus .4em}
$$\sum_{\ell=1}^\infty\int_{E_\ell}F_2\,d\mu=
\int\limits_{\{F_2\ge \frac{\sqrt n}4\}\cap X'}
\!\!\!F_2 \,d\mu
\ge\frac{\sqrt n}4$$
so that at least one of the following cases takes place:\nopagebreak

{\samepage
 \begin{itemize}
\item[\bf (i)]  $\sum\limits_{1\le \ell< n^{\varepsilon_2}}\,
     \int\limits_{E_\ell}F_2\ge {\sqrt n}/8$;
\item[\bf (ii)] $\sum\limits_{\ell\ge n^{\varepsilon_2}}\,
     \int\limits_{E_\ell}F_2> {\sqrt n}/8$.
 \end{itemize}
}

Assume first that {\bf(i)} holds.
Define the following index set
$${\cal L}:=\Big\{1\le
\ell<n^{\varepsilon_2}:
\int_{E_\ell}F_2\ge\frac{n^{\frac12-\varepsilon_2}}{16}\Big\}$$
and notice that
$\sum\limits_{\ell\in{\cal L}}\int_{E_\ell}F_2\ge\frac{\sqrt n}{16}$.
Denote
$$\mu_\ell:=\mu E_\ell;\qquad\qquad \rho_\ell:=\frac{\sqrt n}4 2^{\ell-1}.
$$
It is easy to see that in the case~{\bf(i)}
we have
$$\sum_{\ell\in{\cal L}}\mu_\ell\rho_\ell\ge\frac{\sqrt n}{32};
\qquad\qquad \mu_\ell\rho_\ell\ge\frac{n^{\frac12-\varepsilon_2}}{32}
\qquad\ \text{for all~}\ell\in{\cal L}.$$

Further, on steps 5--11 we shall deal with the case~{\bf(i)} only.
The case~{\bf(ii)} is simpler and we shall consider it on the
final step~12.
\pagebreak[1]



\smallskip
\noindent
{\bf Step 4.}
Assume ${\cal J}_x\subset \{1,\dots,n\}$
is an index set (which may depend on~$x\in X$),
satisfying
${|{\cal J}_x|\ge n-n^{\frac12+2\epsilon}}$ (${\epsilon<1/4}$),
and for each~$k_0$ the set ${\{x\in X: k_0\in{\cal J}_x\}}$
is $\mu$-measurable.
Then for an arbitrary set of signs ${\theta_k=\pm1}$, $k=1,\dots,n$,
it follows that
\begin{align*}
\int_X \big|\sum_{k\in{\cal J}_x}\theta_kf_k(x)\big|d\mu(x)&\le
\int_X\big|\sum_{k=1}^n\theta_kf_k(x)\big|d\mu(x)
+\int_X \big|\!\sum_{k\in{\cal J}^c_x}\theta_kf_k(x)\big|d\mu(x)\\
&\le Mn^{\frac12+p_1}+
\int_X\sqrt{|{\cal J}^c_x|}\cdot\Big(\sum_{k=1}^n|f_k(x)|^2\Big)^{1/2}d\mu(x)\\
&\le Mn^{\frac12+p_1}+Mn^{\frac34+p_2+\epsilon}
\!\le2Mn^{\frac12+\max\{p_1,p_2+\frac14+\epsilon\}}\!.
\end{align*}
Therefore, by Lemma~\ref{geomlem}
for arbitrary coefficients~$\{a_k\}_1^n$ we have
$$\int_X\big|\sum_{k\in{\cal J}_x}a_kf_k\big|d\mu(x)\le
Cn^{\frac14+\frac12\max\{p_1,p_2+\frac14+\epsilon\}}
\Big(\sum_{k=1}^n |a_k|^2\Big)^{\frac12}
$$
with a constant~${C(M)>0}$.
\pagebreak[1]

\smallskip
\noindent
{\bf Step 5.}
Let us choose the constants $\varepsilon_1,\varepsilon_2,\varepsilon_3>0$
such that the following inequalities take place:
\begin{equation}\label{epsdelta}
\begin{array}{rl}
\vphantom{\sum\limits_k}
2\varepsilon_2+\varepsilon_3 &\le \frac14(1-2p_1-4p_2);\\
\frac{\varepsilon_1}2+2\varepsilon_2+\varepsilon_3
&\le\frac32\big(\frac1{12}-p_2\big).
\end{array}
\end{equation}
The reasons for such a choice will be clear soon.
To be definite we could set
$$
2\varepsilon_2=\frac{\varepsilon_1}2=\varepsilon_3:=
\min\big\{\frac12\big(\frac1{12}-p_2\big),\frac18(1-2p_1-4p_2)
\big\}.
$$
The constraints imposed on $p_1$ and~$p_2$
ensure  that $\varepsilon_1$, $\varepsilon_2$ and $\varepsilon_3$
are positive.

Consider some $E_\ell$, $\ell\in{\cal L}$.
For~${\bar{x}=(x_i)_1^m\in (E_\ell)^m}$ set
\begin{align*}
{\cal I}(x_1,\dots,x_m)&:=\bigcap_{j=1}^m{\cal A}_{x_j};\\
{\varphi}(x_1,\dots,x_m)&:=
\frac1{m^2}\sum_{\att{j,k=1}{j\ne k}}^m
\Big|\sum_{i\in {\cal I}(\bar{x})}f_i(x_j)f_i(x_k)\Big|.
\end{align*}
Note that $|{\cal I}(\bar{x})|\ge n-mn^{\frac12+\varepsilon_1}
\ge n-n^{\frac12+2\varepsilon_1}$.
Using the fact that the set ${\cal I}(x_1,\dots,x_m)$
and the function~$\varphi(x_1,\dots,x_m)$
are invariant with respect to permutations of
the coordinates $(x_1,\ldots,x_m)$, we get
\begin{multline*}
\frac1{\mu_\ell^m}
\int\limits_{(E_\ell)^m}{\varphi}(\bar{x})d\mu^m(\bar{x})
=\frac{m(m-1)}{m^2\mu_\ell^m} \int\limits_{(E_\ell)^m}
\Big|\! \sum_{k\in {\cal I}(\bar{x})}\! f_k(x_1)f_k(x_2)\Big|
d\mu^{m}(\bar x)
\le\\
\le\frac1{\mu_\ell^m}\!\!
\int\limits_{(E_\ell)^{m-2}}\!\!\!\Big\{ \int\limits_X\!\int\limits_X
\Big|\! \sum_{k\in {\cal I}(\bar{x})}\! f_k(x_1)f_k(x_2)\Big|
d\mu(x_1)d\mu(x_2)\Big\}d\mu^{m-2}(x_3,\dots,x_m).
\end{multline*}
Applying the arguments of Step~4 for the set
${\cal I}(\bar{x})$ б~$\epsilon=\varepsilon_1$, we get
\begin{multline*}
\int_X\int_X \Big|\sum_{k\in {\cal I}(\bar{x})}
f_k(x_1)f_k(x_2)\Big| d\mu(x_1)d\mu(x_2)\le\\
\le
Cn^{\frac14+\frac12\max\{p_1,p_2+\frac14+\varepsilon_1\}}
\cdot\int_X\Big(\sum_{k=1}^n|f_k(x_1)|^2\Big)^{1/2}d\mu(x_1),
\end{multline*}
which in turn implies
\begin{equation}\label{st5}
\frac1{\mu_\ell^m} \int_{(E_\ell)^m}
{\varphi}(\bar{x})d\mu(x_1)\dots d\mu(x_m) \le
\frac{CM}{\mu_\ell^2}
n^{\frac34+p_2+\frac12\max\{p_1,p_2+\frac14+\varepsilon_1\}}.
\end{equation}

We are going to show that there exists a set~$G_\ell\subset(E_\ell)^m$
such that
${\mu^mG_\ell\ge\frac12\mu_\ell^m}$
and
$${\varphi}(\bar{x})\le C_0\rho_\ell^2n^{-\varepsilon_3}$$
for all~${\bar{x}=(x_i)_1^m\in G_\ell}$
with some constants $\varepsilon_3$, $C_0(p_1,p_2,M)>0$.
Indeed, whenever ${\ell\in{\cal L}}$ we have
$\mu_\ell\rho_\ell\ge2^{-5}n^{\frac12-\varepsilon_2}$
which, combined with~(\ref{st5}), implies
$$
\frac1{\mu_\ell^m} \int_{(E_\ell)^m}{\varphi}(\bar{x})d\mu^m(\bar{x})
\le
CM2^{10}{\rho_\ell^2}n^{2\varepsilon_2-1}
n^{\frac34+p_2+\frac12\max\{p_1,p_2+\frac14+\varepsilon_1\}}
\le CM2^{10}\rho_\ell^2n^{-\varepsilon_3},
$$
where the constants
$\varepsilon_1$, $\varepsilon_2$,~$\varepsilon_3$
chosen positive and satisfying the inequalities~(\ref{epsdelta}).
Thus, by Chebyshev's inequality for the set
$$G_\ell:=\big\{\bar x\in E_\ell:
\varphi(\bar x)\le CM2^{11}\rho_\ell^2n^{-\varepsilon_3}\big\}
,\qquad\ell\in{\cal L},$$
we have
$\mu^m
G_\ell\ge\frac12\mu_\ell^m$.
\pagebreak[1]

\smallskip
\noindent
{\bf Step 6.}
For every $\ell\in{\cal L}$ and $(x_j)_1^m\in G_\ell$
define the random vector ${\eta=(\eta_j)_1^m}$
by
$$
\eta_j(\omega):= |{\cal I}(\bar x)|^{-1/2}\!\!
\sum_{i\in {\cal I}(x_1,\dots,x_m)}\!
\xi_i(\omega)f_i(x_j).
$$
We are going to show that
\vspace{0em plus 0em minus .2em}
\begin{equation}\label{st6}
\sup_{v\in{\Bbb R}^m}{\sf P}
\Big\{\max_{1\le j\le m}
|\eta_j+v_j|\le
\alpha\rho_\ell\sqrt{\frac{\log m}{|{\cal I}(\bar x)|}}\Big\}
\le\frac{K_0}{m^q}
\end{equation}
with some positive constants $\alpha$, $K_0$ and~$q$
which depend only on~$p_1$,~$p_2$,~$M$.
Obviously, inequality~(\ref{st6}) holds if
$$
\min_{\theta_j=\pm1}
{\sf P}\Big(\bigcup_{j=1}^m\Big\{\theta_j\eta_j\ge
\alpha\rho_\ell\Big(\frac{\log m}{|{\cal I}(\bar x)|}\Big)^{1/2}
\Big\}\Big)
\ge 1-\frac{K_0}{m^q} $$
(to see that it suffices to take~${\theta_j=-\text{sign}(v_j)}$).
To prove the latter inequality
we show that
$$
\min_{\theta_j=\pm1}
{\sf P}\Big(\bigcup_{j=1}^m\Big\{\theta_j\eta_j
\ge \alpha\sqrt{d_{j,\ell}\log m}
\Big\}\Big) \ge 1-\frac{K_0}{m^q},
\eqno(\ref{st6}')$$
$$
\text{where}\quad
d_{j,\ell}:={\sf D}(\eta_j)=
\frac1{|{\cal I}(\bar x)|}\sum_{i\in {\cal I}(\bar x)}\!|f_i(x_j)|^2\ge
\frac1{|{\cal I}(\bar x)|}\!\sum_{i\in \Lambda(x_j)}\!\!|f_i(x_j)|^2
\ge\frac{\rho_\ell^2}{|{\cal I}(\bar x)|}.
$$
(While estimating $d_{j,\ell}$ we used the definition of $\Lambda_x$ and the fact that
$|{\cal I}(\bar x)|\ge {n-m\cdot n^{\frac12+2\varepsilon_1}\ge|\Lambda_x|}$.)

For a fixed set of signs~$\{\theta_j\}_1^m$ ($\theta_j=\pm1$)
denote
\begin{align*}
U_j
&=U_j(\alpha):=
\big\{\omega\in\Omega:\theta_j\eta_j\ge \alpha\sqrt{d_{j,\ell}
\log m}\big\};\\
v_{j,k} &:=\frac{\theta_j\theta_k}{|{\cal I}(\bar x)|}
\sum_{i\in {\cal I}(\bar x)}f_i(x_j)f_i(x_k),
\end{align*}
and note that by the definition of the set~$G_\ell$
\begin{equation}\label{st6nv}
\frac1{m^2}\sum_{\att{j,k=1}{j\ne k}}^m|v_{j,k}|\le
C_0\rho_\ell^2n^{-\varepsilon_3}|{\cal I}(\bar x)|^{-1}.
\end{equation}

In what follows we are going to demonstrate that
the random variables~$\eta_j$
are ``almost'' Gaussian and, moreover, ``almost''
parewise Gaussian.
Then we apply Lemma~\ref{lem3tex}
(with parameters~${R=n}$,~$P=m$)
and establish for the events~$U_j$ the following inequality:
\begin{equation}\label{l1star2}
\sum_{j,k=1}^m
{\sf P}\big(U_j(\alpha)\cap U_k(\alpha)\big)
\le (1+L_0 m^{-q}) \Big(\sum_{j=1}^m{\sf P}
\big(U_j(\alpha)\big)\Big)^2
\end{equation}
with some constants $L_0(p_1,p_2,M)$, $q(p_1,p_2)>0$,
this inequality with help of Lemma~\ref{lem1}
implies~(\ref{st6}$'$).
Thus, to prove~(\ref{st6})
it suffices to prove~(\ref{l1star2}).
Steps~7 and~8 are devoted to the proof of~(\ref{l1star2}).
\pagebreak[1]

\smallskip
\noindent
{\bf Step~7.}
Let $(h_j)_1^m$ denote a Gaussian vector with
zero mean and the covariances: ${{\sf E}h_jh_k=v_{j,k}}$ and
${\Psi_j=\Psi_j(\alpha)
:=\{{\omega:h_j>\alpha\sqrt{d_{j,\ell}\log m}}\}}$.
Note that
$$
{{\sf P}(\Psi_j)=(2\pi d_{j,\ell})^{-1/2}
\int_{\alpha\sqrt{d_{j,\ell}\log m}}^\infty
\exp\big(-\frac{y^2}{2d_{j,\ell}}\big)dy}.
$$
We are going to apply Lemma~\ref{lem3tex}
to the vectors~$(\eta_j)_1^m$ and~$(h_j)_1^m$.
In order to check the conditions~(\ref{lem3usl1}),~(\ref{lem3usl2})
we shall use Proposition~\ref{rot}.

Given
fixed $\bar x=(x_j)_1^m\in G_\ell$
and a set of signs~$\{\theta_j\}_1^m$, ${\theta_j=\pm1}$,
let us
apply Proposition~\ref{rot}
to the set of random variables
${\big\{\theta_jf_i(x_j)\xi_i(\omega)\big\}_{i\in {\cal I}(\bar x)}}$
for each ${1\le j\le m}$.
Here ${N=|{\cal I}(\bar x)|}$ and
\begin{align*}
m_3 &=
\frac1{|{\cal I}({\bar x})|}
\sum_{k\in {\cal I}({\bar x})}{\sf E}|\xi_k|^3|f_k(x_j)|^3
\le \frac{M^3}{|{\cal I}(\bar x)|} \sum_{k\in {\cal I}(\bar x)}|f_k(x_j)|^3\\
\lambda &=V=
\frac1{|{\cal I}(\bar x)|}\sum_{k\in {\cal I}({\bar x})} {\sf E} |\xi_k|^2 |f_k(x_j)|^2
 =d_{j,\ell}.
 \end{align*}
Taking into account (\ref{tridva2}), we get
$$
|{\cal I}({\bar x})|^{-1/2}m_3 \lambda^{-3/2}\le 3M^3n^{-\varepsilon_1}
$$
and, consequently, by Proposition~\ref{rot} applied
for the variable~$\theta_j\eta_j$ it follows that
\begin{align}
|{\sf P}(U_j(\alpha))- {\sf P}(\Psi_j(\alpha))|
&\equiv
\big|{\sf P}
(U_j(\alpha))-\Big(\frac{1}{2\pi d_{j,\ell}}\Big)^{1/2}
\hspace{-.1em plus 0em minus .1em}
\int_{\alpha\sqrt{d_{j,\ell}\log m}}^\infty
\hspace{-.2em plus 0em minus .2em}
\exp\Big\{\frac{-y^2}{2d_{j,\ell}}\Big\}dy\big|
\notag\\
\label{l1tex1}
&\le3K_1(1)M^3n^{-\varepsilon_1}.
\end{align}
\pagebreak[1]

\smallskip
\noindent
{\bf Step 8.}
As in the proof of Lemma~\ref{lem2} set
$$\sigma_1:=\big\{(j,k): 1\le j\ne k\le m,
|v_{j,k}|<\frac{\rho_\ell^2}{8|{\cal I}(\bar x)|}\big\}.
$$
By Chebyshev's inequality it follows that
(see~(\ref{st6nv}))
\begin{equation}\label{l1tex3}
|\sigma_1^c|\le
\frac{8|{\cal I}(\bar x)|}{\rho_\ell^2}
\sum_{\att{j,k=1}{j\ne k}}^m|v_{j,k}|
\le {8C_0m^2n^{-\varepsilon_3}}
\le 8C_0m^{2-{\varepsilon_3}/{\varepsilon_1}}.
\end{equation}
(We used the assumption that $m\le n^{\varepsilon_1}$).

To apply Lemma~\ref{lem3tex}
and, thus, to prove~(\ref{l1star2})
it remains to demonstrate that
an estimate of type~(\ref{lem3usl2})
holds for all ${(j,k)\in\sigma_1}$.

For a fixed pare $s=(j,k)\in \sigma_1$
consider the following set of random vectors in~${\Bbb R}^2$:
$$
\Big\{
\Big(\frac{\theta_j\xi_i(\omega)f_i(x_j)}{\sqrt{d_{j,\ell}}},
\frac{\theta_k\xi_i(\omega)f_i(x_k)}{\sqrt{d_{k,\ell}}}
\Big)\Big\}_{i\in {\cal I}(x_1,\dots,x_m)}.
$$
Apply Proposition~\ref{rot} to this set.
Here
\begin{align*}
m_3^s
&=\!\frac1{|{\cal I}(\bar x)|} \sum_{i\in {\cal I}(\bar x)}\!\!
{\sf E}|\xi_i|^3
\Big(\frac{|f_i(x_j)|^2}{d_{j,\ell}}
 +\frac{|f_i(x_k)|^2}{d_{k,\ell}}
\Big)^{3/2};\\
V^s
&=\frac1{|{\cal I}(\bar x)|}
\left(\!
\mbox{$
\begin{array}{l}
\sum\limits_{i\in {\cal I}(\bar x)\vphantom{\frac1I}}
\frac{|f_i(x_j)|^2}{d_{j,\ell}}\\
\frac{\theta_j\theta_k}{\sqrt{d_{j,\ell}d_{k,\ell}}}
\sum\limits_{i\in {\cal I}(\bar x)}\!\!\!f_i\!(x_j)f_i(x_k)
\end{array}
$}
\hspace{-1em}
\mbox{$
\begin{array}{r}
\frac{\theta_j\theta_k}{\sqrt{d_{j,\ell}d_{k,\ell}}}
\sum\limits_{i\in {\cal I}(\bar x)\vphantom{\frac1I}}\!\!\! f_i(x_j)f_i(x_k)\\
\sum\limits_{i\in {\cal I}(\bar x)}
 \frac{|f_i(x_k)|^2}{d_{k,\ell}}
\end{array}
$}
\right).
\end{align*}
Therefore,
$$
V^s=
\left(
\begin{array}{cc}
1 & \frac{v_s}{\sqrt{d_{j,\ell}d_{k,\ell}}} \\
\frac{v_s}{\sqrt{d_{j,\ell}d_{k,\ell}}} & 1
\end{array}\right).
$$
Since
$s\in\sigma_1$,
we can estimate:
$$
\det V^s=1-\frac{v_s^2}{d_{j,\ell}d_{k,\ell}}
\ge 1-\frac{1}{64}
\ge\frac12.
$$
Since the matrix $V^s$
is positive definite
both its eigenvalues: ${\lambda_2\ge\lambda_1}$ are positive.
Therefore, taking into account that
${\lambda_1+\lambda_2}=\text{\rm{trace\ }} V^s=2$,
we get
$$
\frac12\le\det V^s=\lambda_2\lambda_1
\le 2\lambda_1
$$
so that $\lambda_1>1/4$.
Taking into account~(\ref{tridva2})
and the definition of~$d_{j,\ell}$, ${\cal I}(\bar x)$,
we estimate
$$
m^s_3 \le 8M^3|{\cal I}(\bar x)|^{1/2}\cdot 3n^{-\varepsilon_1}.
$$

Thus, by Proposition~\ref{rot} we have
\begin{align}
\Big|&{\sf P}( U_j\cap  U_k)-
\frac1{2\pi\sqrt{\det V^s}}\!
\hspace{-.2em plus 1pt minus .5em}
\int\limits_{\alpha\sqrt{d_{j,\ell}\log m}}^\infty\,
\int\limits_{\alpha\sqrt{d_{k,\ell}\log m}}^\infty
\hspace{-.4em }
e^{-\frac12(Y,(V^s)^{-1}Y)}dy_1dy_2
\Big|\equiv  \notag\\
&\label{l1tex2}
\equiv
\big|{\sf P} ( U_j\cap  U_k)-{\sf P}(\Psi_j\cap \Psi_k)\big|
\le \frac{K_1(2)m_3^s\lambda_1^{-3/2}}{|{\cal I}(\bar x)|^{1/2}}
< \frac{200M^3K_1(2)}{n^{\varepsilon_1}}.
\end{align}

Now, equipped with the estimates~(\ref{l1tex1}), (\ref{l1tex3})
and~(\ref{l1tex2}) we can apply Lemma~\ref{lem3tex}
with parameters~${R=n}$, ${P=m}$,
${\delta_1=\delta_2=\varepsilon_1}$ and
${\delta=\delta_3={\varepsilon_3}/{\varepsilon_1}}$
to derive~(\ref{l1star2})
and, thus, to prove~(\ref{st6})
with some\footnote{In this work we do not try to choose
\samepage
the power parameter~$q$ optimally.
However, note if~${\alpha^2=\frac1{4}(\frac1{12}-p_1)}$
than Lemma~\ref{lem3tex} ensures~(\ref{l1star2}) and~(\ref{st6})
with ${q=\frac14(\frac1{12}-p_1)}$.
}
fixed ${\alpha>0}$,~${q(p_1,p_2)>0}$.
\pagebreak[1]

\smallskip
\noindent
{\bf Step 9.}
The aim of this step is to prove that
\begin{equation}\label{l1st11}
{\sf P}\Big(\bigcup_{j=1}^m
\big\{\omega:\big|\sum_{i=1}^n\xi_i(\omega)f_i(x_j)\big|
\le\alpha\rho_\ell\sqrt{\log m}\big\}\Big)
\le K_0m^{-q}
\end{equation}
for $(x_j)_1^m\in G_\ell$, $\ell\in{\cal L}$.
This inequality will easily follow from~(\ref{st6})
and the following

\begin{lem}\label{tver}
Let $\eta=(\eta_j)_1^m$ and $\eta^c=(\eta^c_j)_1^m$
be independent random vectors and let ${B\subset {\Bbb R}^m}$
be an open or closed set, satisfying
$$\sup_{v\in{\Bbb R}^m}
{\sf P}\big\{\eta+v\in B\big\}\le p$$
with some fixed ${p\in(0,1)}$.
Then
$
{\sf P}\big\{\eta+\eta^c\in B\big\}\le p.
$
\end{lem}

\smallskip
\noindent
{\bf Proof of Lemma \ref{tver}.} Let
$\chi(v,w)$ denote the characteristic function
of the set
${\{(v,w):v+w\in B\}\subset{\Bbb R}^{2m}}$, then
$$
{\sf P}\big\{\eta+\eta^c\in B\big\}
=\int_{{\Bbb R}^m}\int_{{\Bbb R}^m}
\chi(v,w)\, dF_{\eta|\,\eta^c}(v|\,w)dF_{\eta^c}(w),
$$
where $F_{\eta|\,\eta^c}(v|\,w)$ is the
conditional distribution of vector~$\eta$
given~$\eta^c$
(see~\cite{shir} for the definition),
$F_{\eta^c}(w)$ is the distribution of~$\eta^c$.
Since the vectors~$\eta$ and~$\eta^c$
are independent, we have ${F_{\eta |\,\eta^c}(v|\,w)=F_{\eta}(v)}$ and
\begin{align*}
{\sf P}\big\{\eta+\eta^c\in B\big\}
&= \int_{{\Bbb R}^m}\int_{{\Bbb R}^m}
\chi(v,w)\, dF_{\eta}(v)dF_{\eta^c}(w) \\
&=\int_{{\Bbb R}^m} {\sf P}\big\{\eta\in B-w\big\}
\,dF_{\eta^c}(w)\le p
\end{align*}
Lemma~\ref{tver} is proved.$\Box$
\pagebreak[1]

\smallskip
Now, to deduce (\ref{l1st11})\footnote{%
If we were only interested in estimates for
the
\samepage
{\it uniform}  norm of random polynomials
(Corollary~\ref{l1uniform}),
then practically we could  finish the proof on this step.
Indeed, inequality~(\ref{l1st11})  with~$m=[n^{\varepsilon_1}]$
proves Corollary~\ref{l1uniform} for the case~{\bf(i)}
from Step~3. If the case~{\bf(ii)} takes place,
then it suffices to apply Proposition~\ref{rot}
to~$\{f_i(x)\xi_i\}_{i\in {\cal A}(x)}$ at a single point
${x\in E_\ell}$,~${\ell>n^{\varepsilon_1}}$,
and, by Lemma~\ref{tver}, derive the desired estimate
with a ``great reserve.''}
from~(\ref{st6}) it suffices to notice that
the random vectors, defined by
$$\eta_j:=
\frac1{\sqrt{|{\cal I}(\bar x)|}}\sum_{i\in {\cal I}(\bar x)} \xi_if_i(x_j),
\qquad\
\eta^c_j:=\frac1{\sqrt{|{\cal I}(\bar x)|}}\sum_{i=1}^n
\xi_if_i(x_j)-\eta_j,
$$
are independent and apply Lemma~\ref{tver}
to these vectors with parameter ${p=K_0m^{-q}}$ and
${B:=\{(y_j)_1^m\in{\Bbb R}^m: |y_j|\le\alpha\rho_\ell
\sqrt{\frac{\log m}{|{\cal I}(
\bar x)|}}\}}$.
\pagebreak[1]

\smallskip
\noindent
{\bf Step 10.}
For $\ell\in{\cal{L}}$  set
$$
\Omega_\ell:=\Big\{\omega\in\Omega:
\mu^m\big\{\bar x\in G_\ell:
\!\max_{1\le j\le m}\big|
F_n(\omega,x_j)\big|
\ge\alpha\rho_\ell\sqrt{\log m}\big\}
\ge\frac{\mu^mG_\ell}2\Big\}.
$$
From~(\ref{l1st11}) it follows that
$$
\gamma:=
{\sf P}\!\times\!\mu^m
\Big\{(\omega,\bar{x})\in\Omega\times G_\ell:
\!\max_{1\le j\le m}\big|
F_n(\omega,x_j)\big|
\ge\alpha\rho_\ell\sqrt{\log m} \Big\}
\ge\big(1-\frac{K_0}{m^{q}}\big)\mu^mG_\ell.
$$
Thus, we have
$$
{\sf P}(\Omega_\ell)\mu^mG_\ell
+ \big(1-{\sf P}(\Omega_\ell)\big)\frac{\mu^mG_\ell}2
\ge \gamma \ge (1-{K_0}{m^{-q}})\mu^mG_\ell
$$
and, therefore
${\sf P}(\Omega_\ell)\ge1-2K_0m^{-q}$.

We shall need the following
\nopagebreak

\begin{lem}\label{tver2}
Let $T_\ell\ge0$ be numbers, satisfying
$\sum_{\ell=1}^L T_\ell=T$
($L$ may be infinite),
and let $\Omega_\ell$ be events, satisfying ${\sf P}(\Omega_\ell)\ge1-p$.
Then
$${\sf P}\big\{\sum_{\ell=1}^LT_\ell I_{\Omega_\ell}(\omega)
\le\frac T2\big\}\le2p,$$
where $I_{\Omega_\ell}$
are the indicators of $\Omega_\ell$.
\end{lem}

\noindent
{\bf Proof of Lemma \ref{tver2}.}
Set $q={\sf P}(\sum_{\ell=1}^LT_\ell I_{\Omega_\ell}\le\frac T2)$,
then
$$q\frac T2+(1-q)T\ge
{\sf E}\sum_{\ell=1}^L T_\ell I_{\Omega_\ell}\ge T(1-p).$$
Therefore, $q\le2p$.
Lemma \ref{tver2} is proved.$\Box$
\pagebreak[1]

\smallskip
Let us apply Lemma \ref{tver2}
to the numbers $T_\ell:=\mu_\ell\rho_\ell$
and the events~$\Omega_\ell$ for~${\ell\in{\cal L}}$
with the parameter ${p=2K_0m^{-q}}$.
We derive that there exists an event~$\Omega_0$
and a subset ${\cal{L}}_0(\omega)\subset{\cal{L}}$
such that $\Omega_0\subset\bigcap_{\ell\in{\cal{L}}_0}\Omega_\ell$,
$${\sf P}(\Omega_0)\ge1-4K_0m^{-q}$$
and  (see Step~3)
\begin{equation}\label{l1st12}
\sum_{\ell\in{\cal{L}}_0(\omega)}\mu_\ell\rho_\ell
\ge\frac12\sum_{\ell\in{\cal{L}}}\mu_\ell\rho_\ell
\ge\frac{\sqrt n}{64}
\qquad\text{a.s. on}\quad \Omega_0
\end{equation}
\pagebreak[1]

\smallskip
\noindent
{\bf Step 11.}
For  $\ell\in{\cal{L}}$ set
$$\Delta_\ell(\omega):=
\big\{x\in E_\ell:
|F_n(\omega,x)|\ge\alpha\rho_\ell\sqrt{\log m}\big\}.$$
Almost surely on $\Omega_\ell$ we have
\begin{align*}
(\mu E_\ell)^m-\big(\mu E_\ell-\mu\Delta_\ell\big)^m
&=\mu^m\big\{
\bar{x}\in E_\ell^m:\max_{1\le j\le m}
|F_n(\omega,x_j)|\ge\alpha\rho_\ell\sqrt{\log
m} \big\} \\
&\ge\frac{\mu^mG_\ell}2\ge\frac{(\mu E_\ell)^m}4
\equiv\frac{\mu_\ell^m}4.
\end{align*}
Therefore,
$1\ge\frac14+\big(1-\frac{\mu\Delta_\ell}{\mu_\ell}\big)^m$
and
$$
\frac{\mu\Delta_\ell}{\mu_\ell}\ge1-\Big(\frac34\Big)^{1/m}
\ge\frac1{4m}.
$$
There exist some subsets
$\Delta'_\ell(\omega)\subset\Delta_\ell(\omega)$
such that ${\mu\Delta'_\ell=\frac{\mu_\ell}{4m}}$.
Set~${E'=\bigcup_{\ell\in{\cal{L}}_0}\Delta'_\ell(\omega)}$,
${\mu E'=\frac1{4m}\sum_{\ell\in{\cal{L}}_0}\mu_\ell}$.

Now, for almost all $\omega\in\Omega_0$
($\Omega_0$  and ${\cal{L}}_0$ defined on the
previous step)
by Theorem~\ref{kee}
we can estimate the integral-uniform norm
of the random polynomial~$F_n(\omega,x)$
of the type~(\ref{polinom2}),
using its ``relative'' norm~${\|\cdot\|_m^*}$
(see~(\ref{norma3})). We have
\begin{align*}
\|F_n(\omega,\cdot)\|_{m}^*
&\ge
\frac1{\mu E'(\omega)}
\Big(
\int_{E'(\omega)}
|F_n(\omega,x)|d\mu(x)\Big)
\Big(1-\big(1-
\mu E'(\omega)\big)^m\Big)\ge\\
&\ge
\frac1{\mu E'(\omega)}
\Big(\sum_{\ell\in{\cal{L}}_0(\omega)}
\frac{\mu_\ell}{4m}\alpha\rho_\ell\sqrt{\log m}\Big)
\Big(1-\big(1-\mu E'(\omega)\big)^m\Big)=\\
&=
\alpha\sqrt{\log m}\Big(
\sum_{\ell\in{\cal{L}}_0(\omega)}\mu_\ell\rho_\ell\Big)
\cdot t_0^{-1}\Big(1-\big(1-\frac{t_0}{4m}\big)^m\Big),
\end{align*}
where $t_0=\sum_{\ell\in{\cal{L}}_0}\mu_\ell$.
Taking into account~(\ref{l1st12}) and that
${1-(1-\frac t{4m})^m\ge t/8}$
for $t\in(0,1)$ and~$m\ge1$, by Theorem~\ref{kee}  (see~(\ref{kl}))
we get
$$
\|F_n(\omega,\cdot)\|_{m,\infty}\ge\|F_n(\omega,\cdot)\|_{m}^*
\ge 8^{-3}\alpha\sqrt{n\log m}
\qquad\text{a.s. on}~\Omega_0.
$$
This, combined with  ${\sf P}(\Omega_0)\ge1-{4K_0}{m^{-q}}$,
proves~(\ref{l1th1}) if the case{\bf~(i)}
takes place (see Step~3).
\pagebreak[1]

\smallskip
\noindent
{\bf Step 12.}
In order to finish the proof of the theorem,
it remains to consider the case{\bf~(ii)} from Step~3.
Notice, since the inequality~(\ref{tridva2}) holds
a.e.~on~$X'$ (for~${I_x={\cal A}_x}$),
we can apply Proposition~\ref{rot}
(in the one-dimensional case)
for the sum of random variables
${\eta_x:={|{\cal A}_x|}^{-1/2}
\sum_{i\in {\cal A}_x}\xi_if_i(x)}$.
We get
$$
\sup_{v\in{\Bbb R}}
{\sf P}\Big\{|\eta_x+v|\le\gamma\sqrt{d(x)}\Big\}=
(2\pi)^{-1/2}\int_{-\gamma}^{\gamma}
e^{-y^2/2}dy
+O(n^{-\varepsilon_1}),
$$
where $d(x):=|{\cal A}_x|^{-1}\sum_{i\in {\cal A}_x}|f_i(x)|^2$.
Set $\gamma=m^{-\varepsilon_1}$,  then
for almost all~${x\in E_\ell}$, ${\ell\ge1}$, we get
$$
\sup_{v\in{\Bbb R}}
{\sf P}\Big\{|\eta_x+v|\le\sqrt{d(x)}
m^{-\varepsilon_1}\Big\}
\le K'm^{-\varepsilon_1}
$$
with a constant $K'(p_1,p_2,M)>0$.
Taking into account that~${|{\cal A}_x|d(x)\ge\rho_\ell^2}$
for ${x\in E_\ell}$, by Lemma~\ref{tver} we get
$$
{\sf P}\Big\{\big|\sum_{i=1}^n\xi_i f_i(x)\big| \le
\rho_\ell m^{-\varepsilon_1}\Big\}
\le K'm^{-\varepsilon_1}
\qquad
\mbox{for all $\ell\ge1$ a.e. on $E_\ell$}.
$$

For each $\ell\ge n^{\varepsilon_2}$ define the sets:
$$
\Omega'_\ell:=\Big\{\omega\in\Omega:
\mu\big\{x\in E_\ell:
\big| F_n(\omega,x)\big|
\ge \rho_\ell m^{-\varepsilon_1}\big\}
\ge\frac{\mu E_\ell}2\Big\}.
$$
As on Step~10,
it is easy to prove that there exist
events~$\Omega'_0$ and an index set
${{\cal{L}}_0(\omega)\subset\{\ell:\ell\ge n^{\varepsilon_2}\}}$
such that
${\Omega'_0\subset\bigcap_{\ell\in{\cal{L}}_0}\Omega'_\ell}$,
${{\sf P}(\Omega'_0)\ge1-{4K'}{m^{-\varepsilon_1}}}$
and
\begin{equation}\label{l1last}
\sum_{\ell\in{\cal{L}}_0(\omega)}\mu_\ell\rho_\ell
\ge\frac1{2}\sum_{\ell\ge n^{\varepsilon_2}}\mu_\ell\rho_\ell
\ge\frac{\sqrt n}{32}
\end{equation}
a.s. on~$\Omega'_0$. As on Step~11 set
$$
\Delta_\ell(\omega):=\big\{x\in E_\ell:|F_n(\omega,x)|\ge
m^{-\varepsilon_1}\rho_\ell\big\},
\quad\ \ell\ge n^{\varepsilon_2}.
$$
By definition of $\Omega'_\ell$
we have $\mu\Delta_\ell\ge{\mu E_\ell}/2\equiv{\mu_\ell}/2$
a.s.~on~$\Omega'_\ell$,
so there exist subsets
${\Delta'_\ell(\omega)\subset\Delta_\ell(\omega)}$
such that ${\mu\Delta'_\ell=\mu_\ell/2}$.
Set ${E':=\bigcup_{\ell\in{\cal{L}}_0}\Delta'_\ell(\omega)}$.
Estimate the ${\|\cdot\|_m^*}$-norm of random polynomial~(\ref{polinom2})
(see~(\ref{norma3}))
a.s.~on~$\Omega'_0$ as follows
\begin{align*}
\|F_n(\omega,\cdot)\|_{m}^*
&\ge
\frac1{\mu{E'(\omega)}}
\Big(
\int_{E'(\omega)}
|F_n(\omega,x)|d\mu(x)\Big)
\Big(1-\big(1-\mu{E'(\omega)}\big)^m\Big)\\
&\ge
\frac{2}{\sum_{\ell\in{\cal{L}}_0(\omega)}\!\mu_\ell}
\Big(\sum_{\ell\in{\cal{L}}_0(\omega)}
\frac{\mu_\ell}{2}\rho_\ell m^{-\varepsilon_1}\Big)
\Big(1-\Big(1-\!\!\!
\sum_{\ell\in{\cal{L}}_0(\omega)}\frac{\mu_\ell}{2}\Big)^m\Big)\\
&=
m^{-\varepsilon_1}\Big(
\sum_{\ell\in{\cal{L}}_0(\omega)}\mu_\ell\rho_\ell\Big)
\cdot\frac{1-(1-t)^m}{2t},
\end{align*}
where
$t=\sum_{\ell\in{\cal{L}}_0(\omega)}{\mu_\ell}/2
<Mn^{p_2}2^{-n^{\varepsilon_2}}$ (see~(\ref{st3})).
Notice that ${(1-y)^m\le1-\frac m2 y}$ for~${y\in(0,\frac1m)}$,
so provided sufficiently large ${n>n_0(p_1,p_2,M)}$
we have ${t^{-1}(1-(1-t)^m)>m/2}$.
Therefore,
combining~(\ref{l1last}) and Theorem~\ref{kee}, we obtain
$$
\|F_n(\omega,\cdot)\|_{m,\infty}>128^{-1}m^{1-\varepsilon_1}\sqrt{n}
\qquad \ \ \text{a.s. on\ } \Omega'_0.
$$
The proof of Theorem~\ref{l1th}$'$ is completed.\pagebreak[2]

\begin{rem}\label{remgenth2}
(An analog of a remark from \cite{kt1}).
In the statement of Theorem~\ref{l1th}
the condition of uniform boundness of
the third moments ${{\sf E}|\xi_i|^3\le M^3}$
can be relaxed to a weaker condition:
${{\sf E}|\xi_i|^{2+\varepsilon}\le M}$
with some~${\varepsilon>0}$.
In this case the constants in~(\ref{l1th1})
would depend also on~$\varepsilon$.
In order to prove such a statement
it suffices to apply instead of Proposition~\ref{rot}
with a more precise version of
the central limit theorem
(Corollary~18.3 in~\cite{bhat}).
\end{rem}
\pagebreak[1]


\section{Applications and open problems}\label{applic}

{\bf Applications of the Integral-Uniform Norm.}
In~\cite{sem} Montgomery\-{-Smith} and Semenov
in connection with
their research of strictly singular embeddings
of rearrangement invariant spaces in~$L_1[0,1]$
(i.e.~the spaces whose norms are invariant
with respect to measure invariant changes of variable)
put forward a hypothesis which we formulate
in somewhat simplified form:
\begin{conj}\label{consem2}
For an arbitrary set of functions~$\{f_i\}_{i=1}^n$ from~$L_1[0,1]$
such that~${\|f_i\|_1=1}$,
there exist a set of signs~$\{\theta_i\}_{i=1}^n$, ${\theta_i=\pm1}$,
and a constant~${c_0>0}$
such that
\begin{equation}\label{se0}
\big\|\sum_{i=1}^n \theta_i f_i\big\|'_{2^k}\equiv
\sup_{\att{\Delta\subset [0,1]}{\mu\Delta=2^{-k}}}
\Big\{ 2^k
\int_\Delta \big|\sum_{i=1}^n \theta_i f_i(x)\big| d\mu(x)\Big\}
\ge c_0 \sqrt{nk}
\end{equation}
for all~${k=1,\ldots,n}$ (see~(\ref{norma4})).
\end{conj}

Using Theorem~\ref{l1th} we can show
that the assertion of Conjecture~\ref{consem2} is true,
at least if parameter~$k$ in~(\ref{se0})
varies only in
${1,\ldots,[\log n]}$:

\begin{th}\label{thsem}
For an arbitrary set of functions
$\{f_i\}_{i=1}^n\subset L_1[0,1]$
such that ${\|f_i\|_1=1}$,
there exist a sequence of signs
$\{\theta_i\}_{i=1}^n$, $\theta_i=\pm 1$,
and a constant~${c_0>0}$
such that~(\ref{se0}) holds
for all ${k=1,\dots,[\log n]}$.
\end{th}

The proof of Theorem~\ref{thsem}
almost coincide with the proof of Theorem~5 in~\cite{ja3}.
The only difference is that instead of Theorem~A,
formulated in the Introduction, we need to apply
Theorem~\ref{l1th}.

Conjecture~\ref{consem2}
in the general form  stays open.
Note that for ${k\asymp n^\sigma}$
the inequality~(\ref{se0})
cannot be proved by a random signs argument,
since in most cases the order of the uniform norm
(and, thus the integral-uniform one)
of random polynomials~(\ref{polinom2})
is bounded from above by~$\sqrt{n\log n}$,
e.g. this follows  from
the Salem-Zygmund estimate~(\ref{salz})
(see also Th.~B
from Introduction and Th.~4.3 in~\cite{kah}).

\smallskip
Now, let us show that from an
estimate of the integral-uniform norm
it is possible to
get one for the Marcinkiewicz norm.
Let ${\varphi:[0,1]\to[0,1]}$
be an increasing concave function such that
${\varphi(0)=0}$ and ${\varphi(1)=1}$.
Then the Marcinkiewicz space is defined as
a space of functions on
$[0,1]$,  equipped with the norm
(e.g. see~\cite{sem}):
$$ \|f\|_{M(\varphi)}
:=\sup_{0<t<1}\Big\{\frac1{\varphi(t)}\int_0^t f^*(s)\,ds\Big\},
$$
where $f^*$ is the {\it decreasing rearrangement} of~$f$,
defined by
$$
f^*(s):=\inf\Big\{\tau>0: \mu\big\{x\in[0,1]: f(x)\ge\tau
\big\}<s\Big\}.
$$

As a direct corollary of~(\ref{genth2}) or~(\ref{l1th2})
and the equivalence of the norms ${\|\cdot\|'_{m}}$
and ${\|\cdot\|_{m,\infty}}$
(see Th.~\ref{kee})
one can get an estimate for the Marcinkiewicz norm:

\begin{th}
Assume that for a random polynomial~$F_n(\omega,x)$
of type~(\ref{polinom2})
the estimate~(\ref{l1th2})
takes place for all ${m\le n}$.
And let ${\varphi:[0,1]\to[0,1]}$
be an increasing concave  function
such that  $\varphi(0)=0$, $\varphi(1)=1$.
Then the Marcinkiewicz norm of~$F_n(\cdot,x)$
can be estimated as follows:
$$
{\sf E}\|F_n(\omega,x)\|_{M(\varphi)}\ge
A\sqrt{n}\max_{m=2,\dots,n}
\Big\{\frac{\sqrt{\log m}}{m\varphi(1/m)}\Big\}
$$
with a constant~$A>0$.
\end{th}
\pagebreak[1]




\medskip
{\bf On Possible Generalizations of Theorem~\ref{l1th}.}\nopagebreak

\begin{conj}\label{conja1}
The conclusion of Theorem~\ref{l1th}
stays true for random polynomials~(\ref{polinom2})
with respect to functional systems~$\{f_i\}_1^n$
which satisfy condition~{\bf(d)}
with parameter~${p\in[0,1/2)}$
(in Theorem~\ref{l1th} it is assumed that~${p<1/12}$).
\end{conj}

A somewhat weaker form of Conjecture~\ref{conja1}
could be proved
if we could prove a statement of the following type:

\begin{conj}\label{conja3}
{\samepage
For systems of functions~$\{f_i\}_1^n$
which satisfy~{\bf(d)}
with parameter~${p\in[0,1/2)}$
there exist some constants ${\varepsilon_0\in(0,\frac12)}$, ${L>0}$
such that
$$
\int_X
\sum_{1\le k\le n^{1/2+\varepsilon}}f_{k}^*(x)\,d\mu(x)
\le \text{\rm const}\cdot n^{p+L\varepsilon},
$$
provided ${\varepsilon\in(0,\varepsilon_0)}$,
here $f_k^*(x)$ denote
the decreasing ordered values~$|f_j(x)|$
at a fixed point~${x\in X}$.
}
\end{conj}

By the assertion of Conjecture~\ref{conja3}
one could make some estimates in Step~4
in the proof of Theorem~\ref{l1th}$'$
more efficient.
It is ``optimization''
of Step~4
that is most promising
for prospective generalizations
of Theorem~\ref{l1th} in the direction
of Conjecture~\ref{conja1}.

 \end{document}